\documentclass[final]{siamltex}
\usepackage{booktabs}
\usepackage[margin=.9in]{geometry}

\usepackage[mathscr]{euscript}
\usepackage{color,amsmath,amssymb,mathrsfs}
\usepackage{graphicx}
\usepackage{algorithmic,algorithm}
\usepackage{tikz,pgfplots}
\usepackage{upgreek}
\usepackage{showkeys,cite}
\usepackage{multirow}
\usepackage{yfonts}
\usepackage{mathtools}
\usepackage[normalem]{ulem}
\usepackage{latexsym}
\DeclareMathAlphabet{\mathpzc}{OT1}{pzc}{m}{it}
\usepackage{url} 
\usepackage{tikz} 

\usepackage[plainpages=false,colorlinks=true,citecolor=blue,
  filecolor=blue,linkcolor=blue,urlcolor=blue]{hyperref}
  
\newcommand{\email}[1]{\protect\href{mailto:#1}{#1}}

\newcommand{\norm}[2][]{\left\Vert #2\right\Vert_{#1}}

\newcommand{\bs}[1]{\ensuremath{\boldsymbol{#1}}}

\newcommand{\oO}{\tilde{\Omega}}
\newcommand{\os}{\tilde{\sigma}}

\newcommand{\oet}{\tilde{\bs{\eta}}}
\newcommand{\oG}{\tilde{\Gamma}}
\newcommand{\oB}{\tilde{\beta}}
\newcommand{\ou}{\tilde{u}}

\newcommand{\un}{\bs{m}}

\newcommand{\data}{\bs{y}}
\newcommand{\JJ}{\mathcal{J}}

\newcommand{\fwd}{\mathcal{G}}

\newcommand{\GGamma}{\bs{\Gamma}}
\usepackage[utf8]{inputenc}

\begin{document}

\def\addressR{Department of Engineering Science, University of Auckland, Auckland 1010, New Zealand}
\def\addressM{Department of Applied Physics, University of Eastern Finland, Yliopistonranta 1, FIN-70211 Kuopio, Finland}

\author{Ruanui Nicholson\footnotemark[1]
  \and Matti Niskanen\footnotemark[2]}
\renewcommand{\thefootnote}{\fnsymbol{footnote}}
\footnotetext[1]{\addressR\ (\email{ruanui.nicholson@auckland.ac.nz}).}
\footnotetext[2]{\addressM\ (\email{matti.niskanen@uef.fi}).}

\title{Joint Estimation of Robin Coefficient and Domain Boundary for the Poisson Problem}

\maketitle

\begin{abstract}

We consider the problem of simultaneously inferring the heterogeneous coefficient field for a Robin boundary condition on an inaccessible
part of the boundary along with the shape of the boundary for the Poisson problem. Such a problem arises in, for example, corrosion detection, and thermal parameter estimation. We carry out both linearised uncertainty quantification, based on a local Gaussian approximation, and full exploration of the joint posterior using Markov chain Monte Carlo (MCMC) sampling. By exploiting a known invariance property of the Poisson problem, we are able to circumvent the need to re-mesh as the shape of the boundary changes. The linearised uncertainty analysis presented here relies on a local linearisation of the parameter-to-observable map, with respect to both the Robin coefficient and the boundary shape, evaluated at the maximum a posteriori (MAP) estimates. Computation of the MAP estimate is carried out using the Gauss-Newton method. On the other hand, to explore the full joint posterior we use the Metropolis-adjusted Langevin algorithm (MALA), which requires the gradient of the log-posterior. We thus derive both the Fr\'{e}chet derivative of the solution to the Poisson problem with respect to the Robin coefficient and the boundary shape, and the gradient of the log-posterior, which is efficiently computed using the so-called adjoint approach. The performance of the approach is demonstrated via several numerical experiments with simulated data.
\end{abstract}

\section{Introduction}
In this paper, we consider simultaneous estimation of the distributed Robin coefficient field on an inaccessible part of the domain boundary, as well as the shape of the domain for a Poisson problem based on measurements of the potential on an accessible part of the boundary. The problem setup is inspired from physical applications such as corrosion detection or thermal engineering, see for example~\cite{KaupSantosa95,Inglese97,AlessandriniDelPieroRondi03,DivoKassabKapatEtAl05,Shantasiriwan99} among others. There is a considerable amount of literature available, dealing with both theoretical and computational aspects of both the individual problems, i.e.\ estimation of the Robin coefficient or estimation of the domain shape, see e.g~\cite{KaupSantosa95,Inglese97,AlessandriniDelPieroRondi03,Shantasiriwan99,ChaabaneJaoua99,NicholsonPetraKaipio18,FasinoInglese99,Jin07,BaratchartBourgeoisLeblond16,TossavainenVauhkonenHeikkinenEtAl04,AlbuquerqueLaurainSturm20,BanksKojimaWinfree90,FanSunYang08,HettlichRundell96}, as well as of the problem of joint estimation~\cite{CakoniKress07,Bacchelli08,FasinoIngleseMariani07,IngleseMariani04,HarrachMeftahi19,YangSunChu20,BucurGiacomini16,KarageorghisBin-MohsinLesnicEtAl15,Rundell08,CakoniKressSchuft10a,CakoniKressSchuft10b}. We also mention the works \cite{CakoniKress12,CakoniHuKress14}, in which a similar problem is considered using a generalized impedance boundary condition.

We pose the problem in the Bayesian framework \cite{KaipioSomersalo05,CalvettiSomersalo07,Stuart10}, leading to a problem of statistical inference, the solution of which is the posterior density. We initially quantify the posterior uncertainty in the parameters by employing the so-called Laplace approximation to the posterior, see e.g.~\cite{Wong01}, which entails constructing a local Gaussian distribution about the maximum a posteriori (MAP) estimate. To efficiently calculate the MAP estimate we use the Gauss-Newton method~\cite{NocedalWright06}, which necessitates the computation of the Fr\'{e}chet derivative of the potential with respect to both the Robin coefficient and the boundary shape. The first of these derivatives is straightforward~\cite{NicholsonPetraKaipio18}, while on the other hand, determination of the shape derivative is more involved. The specific approach we take to finding the shape derivative is similar to that provided in~\cite{DardeHyvonenSeppanenEtAl13a}, though for an in-depth discussion of shape derivatives, see~\cite{DelfourZolesio11}. We also note the works of \cite{KolehmainenLassasOla05,KolehmainenLassasOlaEtAl13,KolehmainenLassasOla07} for an alternative approach to estimation of the boundary shape for a Poisson problem based on conformal maps.

To accurately characterise the full joint posterior we use a Markov chain Monte Carlo (MCMC) sampling method. Specifically, we employ the Metropolis-adjusted Langevin algorithm (MALA), which requires the gradient of the log-posterior. To this end we utilise the so-called adjoint approach~\cite{BorziSchulz11,HinzePinnauUlbrichEtAl08,Troltzsch10} as an efficient means to compute the gradient.

Current approaches to handling the shape estimation problem can be loosely separated into those which use boundary integral methods (such as the boundary element method), see e.g.~\cite{Rundell08,CakoniKressSchuft10a,CakoniKressSchuft10b}, and those which use volume integral methods (such as the finite element method) ~\cite{BoyleAdlerLionheart12,DardeHyvonenSeppanenEtAl13a,DardeHyvonenSeppanenEtAl13b,HyvonenKarhunenSeppanen10,YangSunChu20}. Both of these approaches, however, suffer from several limitations. First, the boundary-based methods are only useful in a fairly restrictive setting, as the they are typically based on the fundamental solutions. A major drawback of volume-based methods is that they require re-meshing the domain as the estimate of the domain shape changes, or embedding the problem in a larger mesh and using interpolation. The additional computational overheads incurred by having to re-mesh may be negligible when finding the MAP estimate. However, when carrying out MCMC sampling the mesh would need to be adapted for each new proposal, leading to a significant increases in computational time. Furthermore, as alluded to in~\cite{BoyleAdlerLionheart12,DardeHyvonenSeppanenEtAl13b,TossavainenVauhkonenHeikkinenEtAl04}, altering the mesh can place a limit on the feasible boundary movements, and can also introduce spurious artefacts.

Current approaches are further complicated by including (simultaneous) inference for the Robin coefficient over the inaccessible part of the domain. As the problem is posed in the Bayesian framework, we must specify a prior distribution for the Robin coefficient (as well as the boundary shape). However, as the Robin coefficient is defined over the inaccessible part of the boundary, changing the boundary shape would necessitate updating (i.e., recomputing) the prior, adding significant additional cost.

To avoid re-meshing, while also allowing for fairly arbitrary conductivities, we exploit the fact that the Poisson problem is invariant under the {\em push forward}~\cite{Sylvester90,Isakov06,KolehmainenLassasOla05}, as discussed in the following section. This allows for all computations to be carried out in a simple predefined {\em reference domain} with a fixed mesh, though we point out that the resulting forward problem will generally be transformed to the anisotropic version of the Poisson equation even if the actual conductivity in the true domain is isotropic.

The remainder of the paper is organised as follows. In Section 2, we review the forward problem including a key invariance property, and formulate the required shape derivative based on the invariance property. In Section 3 we briefly review the the Bayesian framework for inverse problems, providing the details on the Laplace approximation, as well as on MALA. In Section 4 we consider three numerical examples, based on different physical situations. Lastly, Section 5 lists the concluding remarks.

\section{The Forward Problem}\label{sec: fwdprob}
Let $\Omega\subset\mathbb{R}^2$, be a bounded domain, the specific nature of which is given below. Furthermore, suppose $\sigma\in L^{\infty}(\Omega)$ is a real and scalar-valued with $\sigma\geq c_1>0 $ for all $\bs{x}\in\Omega$. Assume also, that the boundary can be decomposed as $\partial\Omega= \Gamma_{\rm A}\cup \Gamma_{\rm I} \cup \Gamma_{\rm D}$, with $\Gamma_{\rm A}$, $\Gamma_{\rm I}$, and $\Gamma_{\rm D}$ all having positive measure. Finally, let $E=\left\{e_1,e_2,\dots,e_q\right\}$ denote the set of measurement points with $e_i\in\Gamma_{\rm A}$ for $i=1,2,\dots,q$. The forward problem considered is then, find $\left.u_k\right|_{E}$ such that
\begin{equation}
\begin{aligned}\label{eq: fwd1}
    -\nabla\cdot(\sigma\nabla u_k)&=0 & &\text{in }\Omega\\
    \sigma\nabla u_k\cdot \bs{\eta}&= g_k & &\text{on }\Gamma_{\rm A}\\
    \sigma\nabla u_k\cdot \bs{\eta} + \exp\left(\beta\right)u_k&=0 & &\text{on }\Gamma_{\rm I} \\
    u_k&=0 & &\text{on }\Gamma_{\rm D}
\end{aligned}
\end{equation}
for $k=1,2,\dots,\ell$, and $\bs{\eta}$ the outward facing unit normal. In the case of an electrostatic conductor, $u$ represents the electrostatic potential, $\sigma$ the electrical conductivity, $\exp(\beta)$ the boundary admittance\footnote{this representation is used to ensure positivity of the boundary admittance}, and $g_k$ the applied current densities. As discussed in e.g. \cite{CakoniKress07,Bacchelli08},
the application of several current densities helps ensure uniqueness for the solution of the estimation problem in the continuous setting.

We take the domain to be of the form $\Omega=\left\{(x_1,x_2):\, x_1\in(0,L)\, \text{ and }\, 0<x_2<f(x_1)\right\}$
 with  $f\in C^1_{\rm per}([0,L])$ where $C^1_{\rm per}([0,L]):=\left\{f\in C^1([0,L]) \enspace \vert \enspace f(0)=f(L)\right\}$ and $f(s)>0$ for all $s\in[0,L]$. Then the inaccessible (and possibly corroded) part of the domain is $\Gamma_{\rm I}=\left\{(s,f(s)):\, s\in[0,L]\right\}$. That is to say, the inaccessible part of the boundary $\Gamma_{\rm I}$ is the top of the domain and can be represented as a periodic function. On the other hand, the accessible part of the domain $\Gamma_{\rm A}$ is taken to be the bottom of the domain, i.e., $\Gamma_{\rm A}=[0,L]\times\left\{0\right\}$, while $\Gamma_{\rm D}$ makes up the union of the lateral boundaries. See Figure~\ref{fig: schematic} for a schematic of the forward problem. Inline with the literature \cite{FasinoInglese99,FasinoIngleseMariani07,HarrachMeftahi19}, we choose boundary currents
\begin{align}\label{eq: Neumann}
g_k(s,0)=\sin\left(k 2\pi \frac{s}{L}\right) \quad \text{on }\Gamma_{\rm A}, \quad k=1,2,\dots,\ell.
\end{align}

\subsection{Invariance of the Forward Problem}\label{sec: keyprop}
\begin{figure}[]
\centering
\begin{tikzpicture}    
       \node(image) at (-.09,2.3) {
    \includegraphics[scale=0.1305]{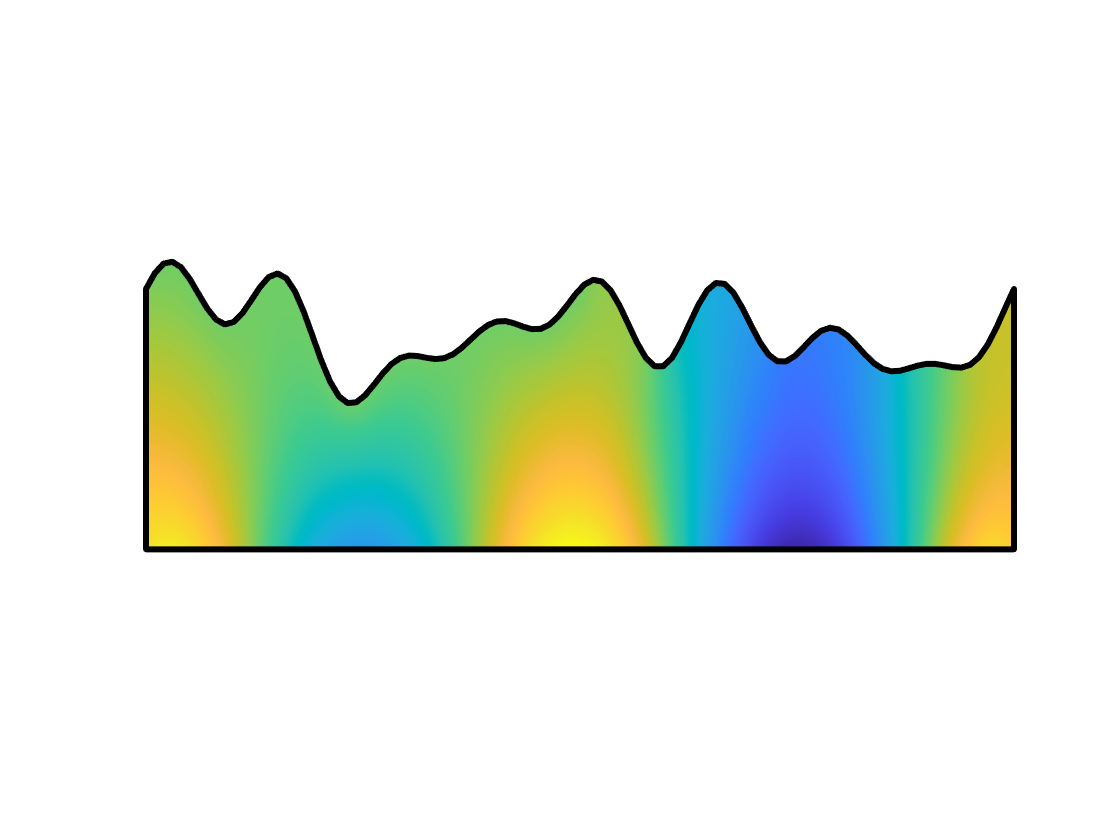}};
    
           \node(image) at (9-.09,2.14) {
    \includegraphics[scale=0.1305]{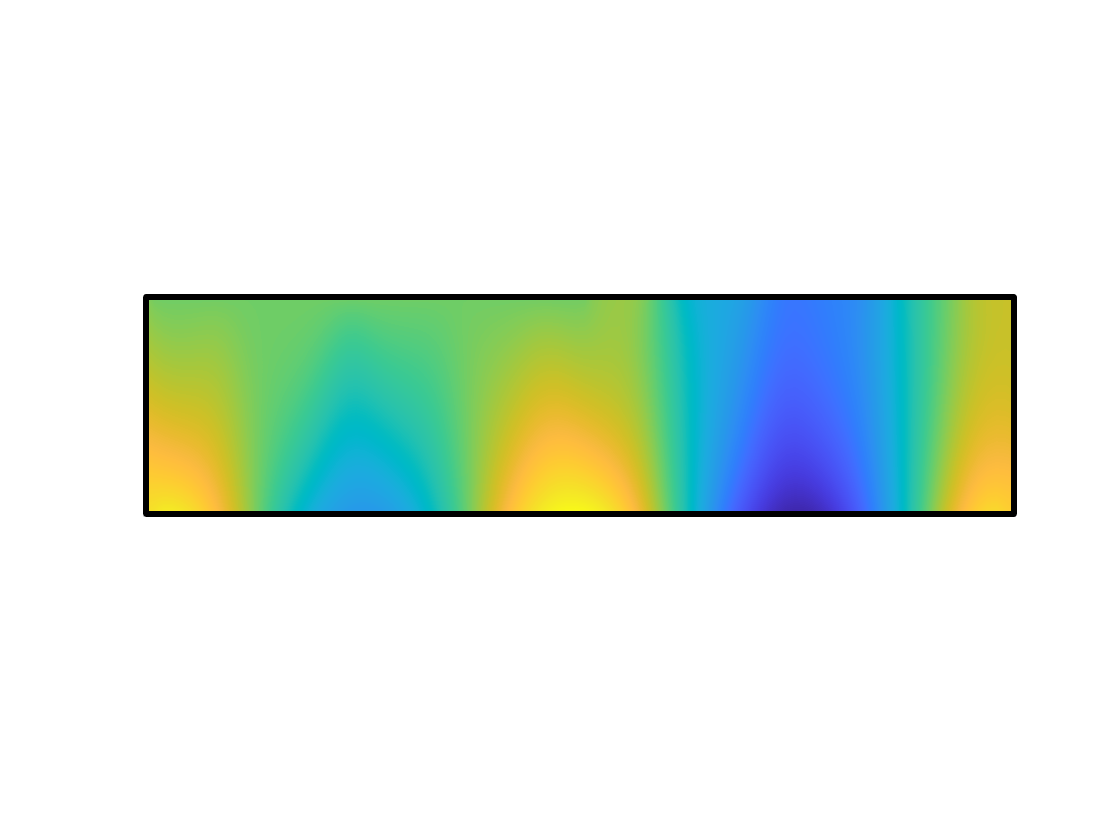}};
    
               \node(image) at (9-.09,-.5) {
    \includegraphics[scale=0.1305]{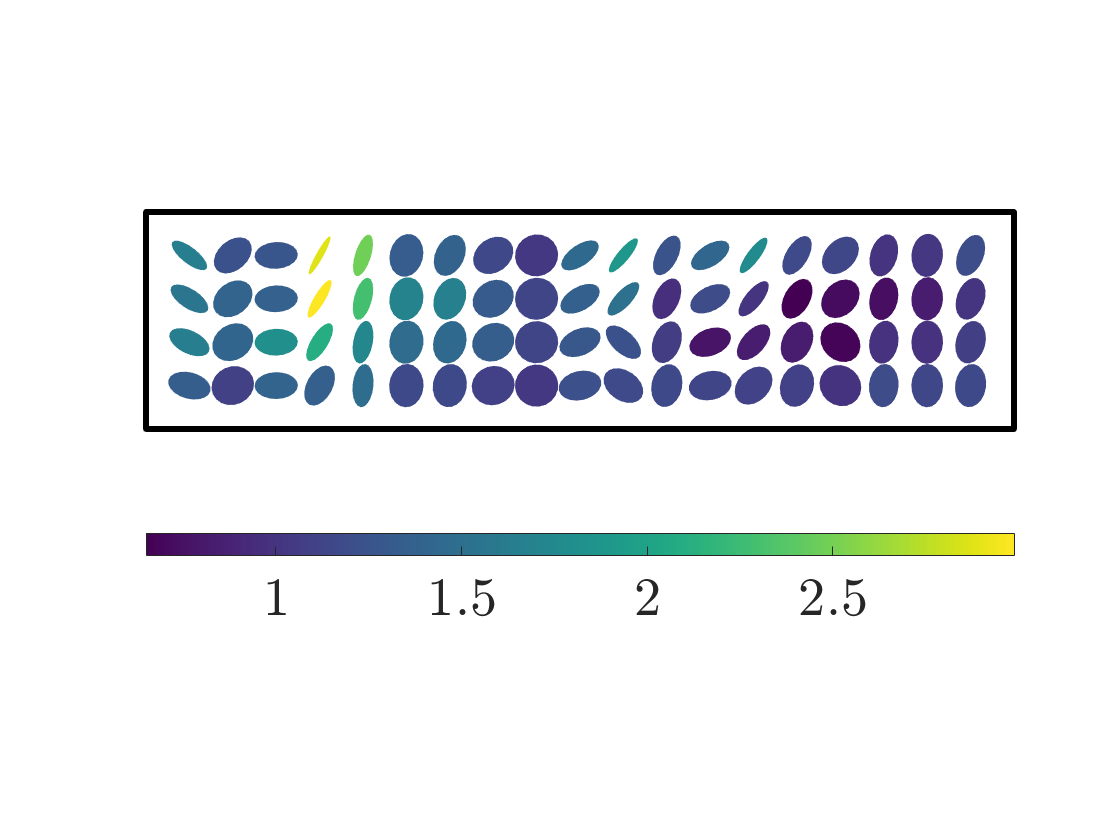}};
           \node(image) at (-.09,-.35) {
    \includegraphics[scale=0.1305]{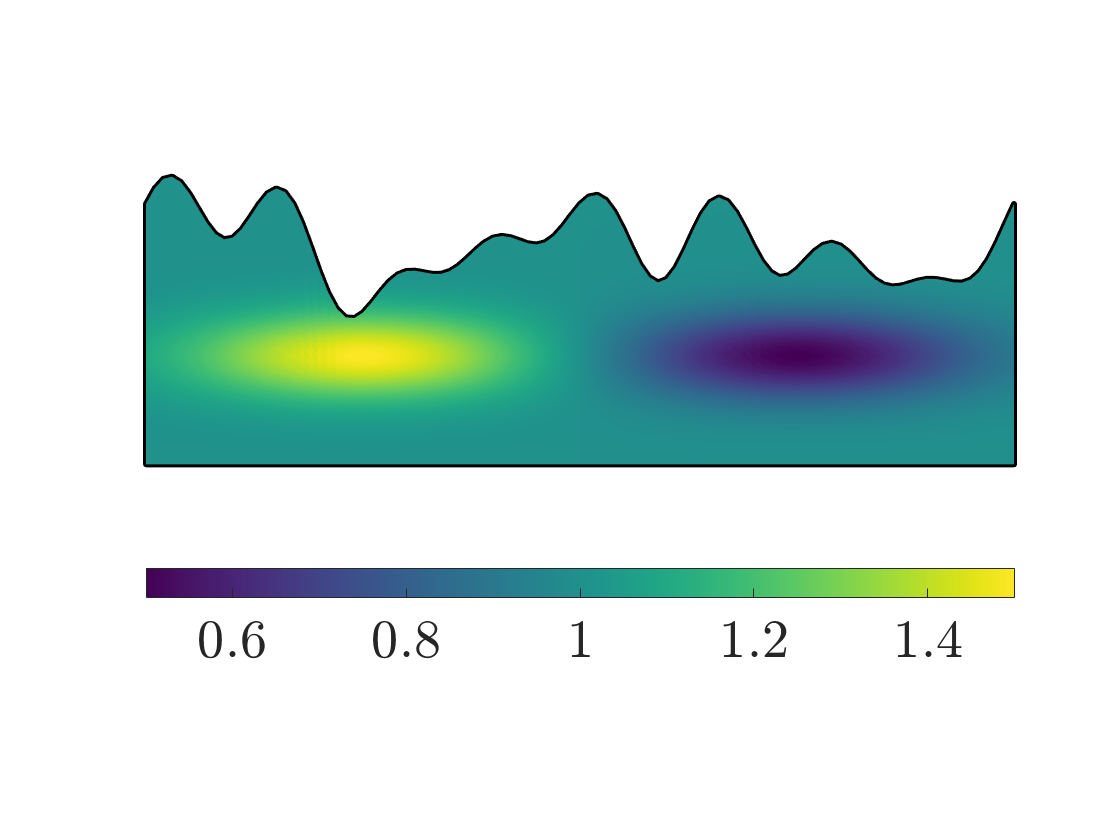}};
    \node[above] at (0,-2) {$\sigma$};
    \node[above] at (9,-2) {$\os$};

    \draw[thick,black] (-2,5.2)--(-2,4)--(2,4)--(2,5.2);
\draw[very thick,black] plot[domain=-2:2,smooth,samples=100] (\x,{5+0.1*sin(2*2/4*pi*(\x+2) r)+0.1*cos(2*2/4*pi*(\x+2) r)+0.1*sin(3*2/4*pi*(\x+2) r)+0.1*cos(6*2/4*pi*(\x+2) r)+0.1*sin(8*2/4*pi*(\x+2) r)});
\draw[thick,black] (8-1,5)--(8-1,4)--(12-1,4)--(12-1,5);
\draw[very thick,black] (8-1,5)--(12-1,5);

        \node[right] at (-1,4.4) {$\Omega$};  
    \node[above] at (0,5.2) {$\Gamma_{\rm I}$}; 
    \node[below] at (0,4) {$\Gamma_{\rm A}$}; 
    \node[right] at (2,4.5) {$\Gamma_{\rm D}$}; 
    \node[left] at (-2,4.5) {$\Gamma_{\rm D}$}; 

        \node[right] at (-1+9,4.4){$\oO$};  
    \node[above] at (0+9,5) {$\oG_{\rm I}$}; 
    \node[below] at (0+9,4) {$\Gamma_{\rm A}$}; 
    \node[right] at (2+9,4.5) {$\oG_{\rm D}$}; 
    \node[left] at (-2+9,4.5) {$\oG_{\rm D}$}; 
    
            \draw[->,thick] (4,2.25)--(5,2.25);
        \node[left] at (4,2.25) {$u$};
        \node[right] at (5,2.25) {$\ou$};
    
        \draw[->,thick] (4,-.25)--(5,-.25); 
        \node[left] at (4,-.25) {$\sigma$};
        \node[right] at (5,-.25) {$\os$};

        \draw[->,thick] (4,4.5)--(5,4.5);  
        \node[left] at (4,4.5) {$\Omega$};
        \node[right] at (5,4.5) {$\oO$};
     \draw[->,thick] (2.5,6.25)--(6.5,6.25);
        \node[above] at (4.5,6.25) {$\psi$};

 \end{tikzpicture}

\vspace{-1cm}

\caption{Illustration of the invariance property of the Poisson problem.  The diffeomorphism, $\psi$, maps the true (corroded) domain, $\Omega$, to a reference (uncorroded) domain $\oO$, which in turn leads to an updated potential, $\ou(\tilde{\bs{x}})=u\circ \psi^{-1}(\tilde{\bs{x}})$, and an updated conductivity $\os$, and boundary conditions (see (\ref{eq: fwd2})-(\ref{eq: bupdates})). The updated conductivity is in general ansitropic. To illustrate this we have plotted the anisotropy of the conductivity using ellipses, the principal axis of which shows the anisotropy while the colour represents the magnitude.}
\label{fig: schematic}
\end{figure}

In the current paper, we denote by $\oO$ the (known) uncorroded reference geometry, which is assumed to be a slab, i.e.,  $\oO=(0,L)\times(0,H)$, with $L\gg H>0$. The Poisson equation, as outlined in (\ref{eq: fwd1}), is known to be invariant under certain deformations of the domain $\Omega$ \cite{Sylvester90,Isakov06,KolehmainenLassasOla05}. More specifically, we let $\psi$ denote the diffeomorphism mapping $\Omega$ to $\oO$ given by,
\begin{align}\label{eq: diffeo1}
    \psi:(x_1,x_2)\mapsto\left(x_1,x_2(f(x_1))^{-1}\right), \quad (x_1,x_2)\in\Omega.
\end{align}
We use $\bs{x}$ to denote points in $\Omega$, while points in $\oO$ will be denoted by $\tilde{\bs{x}}$, i.e., $\bs{x}\in\Omega$ and $\tilde{\bs{x}}\in\oO$. Then, taking into account that such a diffeomorphisms leaves $\Gamma_{\rm A}$ unchanged\footnote{The shape of the accessible part of the domain is assumed known and fixed.}, i.e., $\left.\psi\right|_{\Gamma_{\rm A}}={\rm Id}$ (the identity mapping), the functions $\ou_k(\tilde{\bs{x}})=u_k\circ \psi^{-1}(\tilde{\bs{x}})$ are known to satisfy a transformed Poisson equation \cite{Sylvester90,Isakov06,KolehmainenLassasOla05}:
\begin{equation}
\begin{aligned}\label{eq: fwd2}
    -\nabla\cdot(\os\nabla \ou_k)&=0 & &\text{in }\oO\\
    \os\nabla \ou_k\cdot \oet&=g_k& &\text{on }\oG_{\rm A}\\
    \os\nabla \ou_k\cdot \oet + \exp(\oB)\ou_k&=0 & &\text{on }\oG_{\rm I} \\
    \ou_k&=0 & &\text{on }\oG_{\rm D},
\end{aligned}
\end{equation}
for $k=1,2,\dots,\ell$, where $\oet$ is the outward facing normal of $\oO$, the push forward of $\sigma$ (by $\psi$) is given by
\begin{align}\label{eq: updates}
\os(\tilde{\bs{x}})=\sigma(\bs{x}) \frac{{\bs J}(\bs{x}){\bs J}(\bs{x})^T}{\vert \det({\bs J}(\bs{x}))\vert},\quad {\bs{x}=\psi^{-1}(\tilde{\bs{x}})}
\end{align}
with ${\bs J}$ the Jacobian matrix of $\psi$, and 
\begin{align}\label{eq: bupdates}
\exp\left(\oB(\tilde{s})\right)=\exp\left(\beta(s)\right)\sqrt{1+\left(\frac{df}{ds}\right)^2H^2},\quad \quad {s=\varphi^{-1}(\tilde{s})},
\end{align} 
where $\varphi$ is the restriction of $\psi$ to $\Gamma_{\rm I}$.
It's worth pointing out that even if the conductivity $\sigma$ in $\Omega$ is isotropic, the pushed forward conductivity $\os$ in $\oO$ will in general be anisotropic. We show a schematic
of the invariance property of the Poisson problem in Figure~\ref{fig: schematic}. 

The invariance of the Poisson equation, as outlined above, means that for any domain, $\Omega$, which is diffeomorphic to $\oO$, the solution to (\ref{eq: fwd1}) can be found by solving (\ref{eq: fwd2}). If we only require the solution to (\ref{eq: fwd1}) in a single domain, $\Omega$, there is little use in solving (\ref{eq: fwd2}) instead. However, when considering an estimation problem, such as that considered in the current paper, the forward problem will need to be solved repeatedly, see Section \ref{sec: inf}. Then, since the domain boundary is to be estimated, we will require the solution to (\ref{eq: fwd1}) in a number of different domains. Computing each of these solutions would require re-meshing of the domain, or at least of the inaccessible part of the boundary, $\Gamma_{\rm I}$. This leads to an increased computational overhead, which can become significant if the number of forward solutions required is large. Furthermore, as the Robin boundary coefficient is defined on $\Gamma_{\rm I}$, special care would need to be taken to ensure that spatial prior information, such as smoothness of the coefficient, was kept consistent across each of the domains. This would further increase the computational requirements, and could also lead to convergence issues for the estimates, see e.g. \cite{LassasSaksmanSiltanen09,Stuart10,Bui-ThanhGhattasMartinEtAl13}. To reduce these issues, in the current paper, we propose only ever solving the forward problem in a fixed reference domain, $\oO$, by exploiting the discussed invariance property.

Key to further reducing the computational costs of solving the estimation problem is efficient calculation of derivatives of the forward model with respect to both the shape of the inaccessible part of the domain boundary, $\Gamma_{\rm A}$, and the Robin coefficient, $\beta$, which we now outline.

\subsection{Derivatives}\label{sec: Der1}
In this section, we introduce the generalised derivatives of the forward problem
with respect to both perturbations of the inaccessible part of the boundary $\Gamma_{\rm I}$ and perturbations of the Robin coefficient $\beta$. We point out that a related approach is provided in \cite{DardeHyvonenSeppanenEtAl13a,DardeHyvonenSeppanenEtAl13b}, though the procedure used there requires, among other things, a dual procedure to cope with distributional boundary conditions. In what follows, we will slightly abuse notation and set $g=g_k$ and $u=u_k$ for some fixed $k\in\left\{1,2,\dots,\ell\right\}$. Furthermore, it will be useful to consider writing the weak form of (\ref{eq: fwd1}) as $a(u,v)=F(v)$, with
\begin{align}
a(u,v):=\int_{\Omega}\sigma\nabla u\cdot\nabla v\;d\bs{x}+\int_{\Gamma_{\rm I}}\exp(\beta) u v\;d\bs{s}_{\rm I},\quad\text{and}\quad F(v):=\int_{\Gamma_{\rm A}}g v\;d\bs{s}_{\rm A},
\end{align} 
for $v\in H^1_D(\Omega)$, with $H^1_D(\Omega)=\left\{v\in H^1(\Omega)\vert\enspace \left.v\right|_{\Gamma_{\rm D}}=0\right\}$. 

To calculate the generalised derivative for the boundary we first take $h\in C^1_{\rm per}([0,L])$ and introduce the {\em perturbed problem}: $a^h(u^h,v^h)=F(v^h)$, for which the weak solution, $u^h\in H^1(\Omega^h)$, satisfies
\begin{align}\label{eq: pertprob}
    \int_{\Omega^h}\sigma\nabla u^h\cdot\nabla v^h\;d\bs{x}+\int_{\Gamma^h _{\rm I}}\exp(\beta) u^h v^h\;d\bs{s}_{\rm I}=\int_{\Gamma_{\rm A}}g v^h\;d\bs{s}_{\rm A}\qquad\forall v^h\in H^1_D(\Omega^h),
\end{align}
where $\Omega^h=\left\{(x_1,x_2):\, x_1\in(0,L)\, \text{ and }\, 0<x_2<f(x_1)+h(x_1)\right\}$. We next introduce the diffeomorphism which maps $\Omega^h$ to $\tilde{\Omega}$,  
\begin{align}
\psi^h:(x_1,x_2+h(x_1))\mapsto\left(x_1,(x_2+h(x_1))(f(x_1)+h(x_1))^{-1}\right).
\end{align}

The Gateaux derivative, see e.g. \cite{DelfourZolesio01}, of $u\in H^1(\Omega)$ in the direction $h$ is then
\begin{align}
   \mathsf{D}_hu[\Omega]=\lim_{t\rightarrow 0}\frac{u^{th}\circ\psi^{th}-u\circ\psi}{t},
\end{align}
where the superscript $th$ denotes replacing $h$ with $th$. In what follows, we abuse notation and use the superscript $h$ to denote $th$.

As the potential $u$ can be written as $u=T^{-1}F$ with $T$ the operator associated with a symmetric form $a(u,v)$ \cite{Kato13}, straightforward differentiation yields
\begin{align}\label{eq: spacederiv}
     \mathsf{D}_hu[\Omega]=-T^{-1}\left(\mathsf{D}_h a[\Omega]\right)u+T^{-1} \left(\mathsf{D}_hF[\Omega]\right).
\end{align}
Importantly, due to the nature of the admissible set of boundary perturbations, it follows that $\mathsf{D}_hF[\Omega]=0$. Thus, finding $\mathsf{D}_hL[\Omega]$ renders calculating $\mathsf{D}_hu[\Omega]$ essentially trivial. 

With the two diffeomorphisms, $\psi$ (see (\ref{eq: diffeo1})) and $\psi^{h}$, in hand, we can rewrite both $a$ and $a^{h}$ using the invariance property of the Poisson problem (see Section \ref{sec: keyprop}) as
\begin{align}
    \tilde{a}(\tilde{u},\tilde{v})&=\int_{\oO}\os\nabla\tilde{u} \cdot\nabla \tilde{v}\;d\tilde{\bs{x}}+\int_{\oG _{\rm I}}\exp{(\oB)} \tilde{u}\tilde{v}\;d\tilde{\bs{s}}_{\rm I},\\    \tilde{a}^{h}(\tilde{u}^{h},v^{h})&=\int_{\oO}\os^{h}\nabla\tilde{u}^{h} \cdot\nabla \tilde{v}^{h}\;d\tilde{\bs{x}}+\int_{\oG _{\rm I}}\exp{(\oB^{h})}\tilde{u}^{h} \tilde{v}^{h}\;d\tilde{\bs{s}}_{\rm I},
\end{align}
respectively, where $\os$ and $\exp{(\oB)}$ are defined in (\ref{eq: updates}) and (\ref{eq: bupdates}) respectively, $\ou=u\circ\psi$, while $\ou^{h}=u^{h}\circ \left(\psi^{h}\right)^{-1}$, and 
\begin{align}
   \os^h(\tilde{\bs{x}})=\sigma(\bs{x}^h)\frac{{\bs J}^h(\bs{x}^h) {\bs J}^h(\bs{x}^h)^T}{\vert \det({\bs J}^h(\bs{x}^h))\vert},\quad  \exp\left(\oB^h(\tilde{s})\right)=\exp\left(\beta(s)\right)\sqrt{1+\left(\frac{d\left(f(s)+th(s)\right)}{ds}\right)^2H^2},
\end{align}
with ${\bs J}^h$ the Jacobian matrix of $\psi^h$ and where $\bs{x}^h=\left(\psi^h\right)^{-1}(\tilde{\bs{x}})$ and $s=\varphi^{-1}(\tilde{s})$.
We can then write the derivative as,
\begin{align}
(\mathsf{D}_h\tilde{a}[\Omega])\tilde{u}&=\lim_{t\rightarrow0}\left(\int_{\oO}\os^h\nabla\tilde{u} \cdot\nabla \tilde{v}\;d\tilde{\bs{x}}+\int_{\oG _{\rm I}}\exp(\oB^h) \tilde{u} \tilde{v}\;d\tilde{\bs{s}}_{\rm I}-\int_{\oO}\os\nabla\tilde{u} \cdot\nabla \tilde{v}\;d\tilde{\bs{x}}-\int_{\oG _{\rm I}}\exp(\oB) \tilde{u} \tilde{v}\;d\tilde{\bs{s}}_{\rm I}\right)t^{-1}\\
&=\lim_{t\rightarrow0}\left(\int_{\oO}\left(\os^h-\os\right)\nabla\tilde{u} \cdot\nabla \tilde{v}\;d\tilde{\bs{x}}\right)t^{-1}+\lim_{t\rightarrow0}\left(\int_{\oG _{\rm I}}\left(\exp(\oB^h)-\exp(\oB)\right)\tilde{u} \tilde{v}\;d\tilde{\bs{s}}_{\rm I}\right)t^{-1}\\
&=(\mathsf{D}_{\mathsf{D}_h\tilde{\sigma}[\Omega]}\tilde{a}[\tilde{\sigma}])\tilde{u}+(\mathsf{D}_{\mathsf{D}_hB[\Omega]}\tilde{a}[B])\tilde{u}\label{eq: dad},
\end{align} 
where $\mathsf{D}_h\tilde{\sigma}[\Omega]\in S^2\left(C\left([0,L]\right)\right)$ (i.e., the space of symmetric matrices with elements in $C\left([0,L]\right)$), and where $B=\exp(\oB)$, and $\mathsf{D}_hB[\Omega]\in C\left([0,L]\right).$

In the current set up, due to the assumptions on the domain shape perturbations, each of the derivatives in (\ref{eq: dad}) are straightforward to compute as follows. First,  
since the push forward of $\sigma$ is given by
\begin{align}
\os(\tilde{\bs{x}})=\sigma(\bs{x}) \frac{{\bs J}(\bs{x}){\bs J}(\bs{x})^T}{\vert \det({\bs J}(\bs{x}))\vert}= \sigma(\bs{x})\begin{pmatrix} 
\begin{array}{c|c}
f(x_1) & -\frac{x_2}{f(x_1)}\frac{df(x_1)}{dx_1} \\
\hline
-\frac{x_2}{f(x_1)}\frac{df(x_1)}{dx_1} & \frac{1}{f(x_1)}+\frac{x_2^2}{f^3(x_1)}\left(\frac{df(x_1)}{dx_1}\right)^2
\end{array}
\end{pmatrix},
\end{align}
we have $\mathsf{D}_h\tilde{\sigma}[\Omega]=\sigma(\bs{x})\bs{\Sigma}(\bs{x})h(x_1)$, with 
\begin{align}
\bs{\Sigma}(\bs{x})=\begin{pmatrix}
\begin{array}{c|c}
1 & \left(\frac{x_2}{f^2(x_1)}-\frac{1}{f(x_1)}\right)\frac{df(x_1)}{dx_1}-\frac{x_2}{f(x_1)}\frac{d}{dx_1} \\
\hline
\left(\frac{x_2}{f^2(x_1)}-\frac{1}{f(x_1)}\right)\frac{df(x_1)}{dx_1}-\frac{x_2}{f(x_1)}\frac{d}{dx_1} &\left(\frac{2x_2}{f^3(x_1)}-\frac{3x_2^2}{f^4(x_1)}\right)\left(\frac{df(x_1)}{dx_1}\right)^2+\frac{2x_2^2}{f^3(x_1)}\frac{df(x_1)}{dx_1}\frac{d}{dx_1}-\frac{1}{f^2(x_1)}
\end{array}
\end{pmatrix}.
\end{align}
On the other hand, inline with e.g.~\cite{LechleiterRieder08},  and setting $\hat{h}=\mathsf{D}_h\tilde{\sigma}$, we have
\begin{align}
(\mathsf{D}_{\hat{h}}\tilde{a}[\tilde{\sigma}])\tilde{u}=\int_{\oO}\hat{h}\nabla\tilde{u} \cdot\nabla \tilde{v}\;d\tilde{\bs{x}}.
\end{align}
Next, 
\begin{align}
{\mathsf{D}_hB[\Omega]}=\exp(\beta(s))\left(\frac{H^2\frac{d f(s)}{ds}\frac{d }{ds}}{\sqrt{1+\left(\frac{df(s)}{ds}\right)^2H^2}}\right)h(s),
\end{align}
while finally, by setting $\check{h}={\mathsf{D}_hB[\Omega]}$, we have 
\begin{align}
(\mathsf{D}_{\check{h}}\tilde{a}[B])\tilde{u}=\int_{\oG _{\rm I}}\check{h}\tilde{u} \tilde{v}\;d\tilde{\bs{s}}_{\rm I}.
\end{align}

Notice that in the process of determining the derivatives of the forward problem
with respect to perturbations of the inaccessible part of the boundary, we have also determined the derivatives of the forward problem
with respect to perturbations of the Robin coefficient.  

\section{Inference}\label{sec: inf}
We consider the problem of finding the inaccessible part of the domain and the Robin coefficient as a problem of statistical inference within the Bayesian framework. The solution to the inference problem is the posterior density, i.e., the conditional probability distribution of the parameters conditioned on the measured data. We initially treat the inaccessible part of the domain boundary $\Gamma_{\rm I}$ as a function of $x_1$, and the Robin coefficient $\beta$ as a function defined on $\Gamma_{\rm I}$. The posterior is then expressed using the infinite dimensional version of Bayes' theorem,
\begin{align}\label{eq: rn}
    \frac{1}{Z}\frac{d\mu_{\rm post}}{d\mu_{\rm prior}}\propto \pi_{\rm like}(\data\vert\Gamma_{\rm I},\beta),
\end{align}
where $Z$ is a normalization constant, $\frac{d\mu_{\rm post}}{d\mu_{\rm prior}}$ denotes the Radon-Nikodym derivative of the posterior measure $\mu_{\rm post}$ with respect to the prior measure $\mu_{\rm prior}$, and $\pi_{\rm like}$ denotes likelihood. For an in-depth discussion of (\ref{eq: rn}), and on when the posterior is well-defined see \cite{Stuart10}. 

Although (\ref{eq: rn}) is valid in infinite (and finite) dimensions, we follow the more intuitive form of Bayes’ formula that uses Lebesgue measures and
thus only holds in finite dimensions, see e.g. \cite{Stuart10,Bui-ThanhGhattasMartinEtAl13}. As such, in the following we denote by $\pi_{\rm prior}$ and $\pi_{\rm post}$ the finite dimensional prior and posterior densities (with respect to Lebesgue measure), respectively, and work with the more familiar finite-dimensional version of Bayes' theorem.
In the next section we provide details on the discretisations and parameterisations used for the unknowns, as well as the specific prior distributions used.

\subsection{Parameterisations and the Prior} \label{sec:parametr-prior}
As stated previously, the inaccessible part of the domain boundary is treated as a function of $x_1$ only. Similarly to \cite{DardeHyvonenSeppanenEtAl13a,DardeHyvonenSeppanenEtAl13b}, we take $\Gamma_{\rm I}=\left\{(s,f(s)):\, s\in[0,L]\right\}$ with 
$f\in C^\infty_{\rm per}([0,L])$ and take
\begin{align} \label{eq:Fourier-representation}
    \Gamma_{\rm I}={\color{red}}H\left(1+\alpha_0+\sum_{n=1}^p \alpha_{2n} \cos\left(n2\pi \frac{x_1}{L}\right)+\alpha_{2n-1} \sin\left(n2\pi \frac{x_1}{L}\right)\right), \quad \alpha_0,\alpha_1,\dots.\alpha_{2p}\in\mathbb{R}
\end{align}
As such, the problem of inferring the shape of the inaccessible part of the domain is posed as inferring the Fourier coefficients\footnote{Various other representations, such as wavelets \cite{Daubechies92}, could also be used for the boundary shape.},  $\bs{\alpha}=[\alpha_0,\alpha_1,\dots,\alpha_{2p}]^T\in\mathbb{R}^{2p+1}$.
As in \cite{DardeHyvonenSeppanenEtAl13a,DardeHyvonenSeppanenEtAl13b}, we encode our prior beliefs about the geometry of the domain through a Gaussian prior on $\bs{\alpha}$, with mean $\bs{\alpha}_\ast$ and covariance matrix $\GGamma_\alpha$. As such, the prior density for $\bs{\alpha}$ can be written as 
\begin{align}
    \pi(\bs{\alpha})\propto\exp\left\{-\frac{1}{2}\norm[\GGamma_\alpha^{-1}]{\bs{\alpha}-\bs{\alpha}_\ast}^2\right\}.
\end{align}
Specification of $\bs{\alpha}_\ast$ and $\GGamma_\alpha$ is left to Section \ref{sec: numex}.

We postulate a Gaussian prior measure on the Robin coefficient, i.e., $\mu_\beta=\mathcal{N}(\beta_\ast,\mathcal{C}_\beta)$. Inline with e.g. \cite{FlathWilcoxAkcelikEtAl11,Bui-ThanhGhattasMartinEtAl13,
  PetraMartinStadlerEtAl14}, to ensure the inverse problem is well-posed in infinite dimensions \cite{Stuart10}, while simultaneously promoting smoothness, we use a squared inverse elliptic operator to define the
prior covariance
operator. Specifically, we take $\mathcal{C}_\beta=\mathcal{A}^{-1}$, where $\mathcal{A}$ is the second order elliptic
differential operator defined by
\begin{align}\label{eq: WeakPrior}
\mathcal{A}\beta := -\nabla \cdot (\delta_\beta^{-2} \nabla\beta) +
l^2\delta_\beta^{-2} \beta\quad \text{on } \Gamma_{\rm I},
\end{align}
where $\delta_\beta$ controls the marginal variance and $l$ is the inverse of the correlation length. As discussed in
\cite{KhristenkoScarabosioSwierczynskiEtAl19,DaonStadler18,RoininenHuttunenLasanen14},
suitable boundary conditions should also be applied to reduce so-called boundary
artefacts. To this end, letting $\bs{\eta}_{\rm I}$ denote the outward facing unit vector normal to $\partial\Gamma_{\rm I}$, we equip $\mathcal{A}$ with a Robin boundary condition,
\begin{align}
\delta_\beta^{-2}\nabla\beta\cdot \bs{\eta}_{\rm I}+l\delta_\beta^{-2} \beta=0\quad \text{at } \partial\Gamma_{\rm I},
\end{align}
which leads to a constant (i.e., homogeneous) marginal variance for $\beta$ over $\Gamma_{\rm I}$~\cite{DaonStadler18,RoininenHuttunenLasanen14}.

The Robin coefficient is discretised using continuous piece-wise linear Lagrange basis functions $\left\{\phi_i\right\}_{i=1}^q$, leading to a discrete approximation for $\beta$ of the form $\beta_h=\sum_{i=1}^q\beta_i\phi_i$. Consequentially, the parameters of interest for the Robin coefficient are $\bs{\beta}=[\beta_1,\beta_2,\dots,\beta_{q}]^T\in\mathbb{R}^{q}$, with the resulting discrete representation of the prior covariance
operator, denoted $\GGamma_\beta$, being
\begin{align}
\left[\GGamma_\beta^{-1}\right]_{ij}=\int_{\Gamma_{\rm I}}\phi_i(s) \mathcal{A}\phi_j(s)\; ds\quad i,j\in\{1,2,\dots,q\},
\end{align}
see for example~\cite{Bui-ThanhGhattasMartinEtAl13,PetraMartinStadlerEtAl14,VillaPetraGhattas20}. That is, the covariance matrix can be written as
\begin{align} \label{eq:beta_covariance}
    \GGamma_\beta=\delta_\beta^2(\bs{K}+\bs{M}+\bs{R})^{-1}
\end{align}
where
\begin{align}
    \bs{K}_{ij}=\int_{\Gamma_{\rm I}}\nabla\phi_i\cdot \nabla\phi_j\;ds, \quad \bs{M}_{ij}=\int_{\Gamma_{\rm I}}l^2\phi_i\phi_j\;ds,\quad \bs{R}_{ij}=\int_{\partial\Gamma_{\rm I}}l\phi_i\phi_j\;dt,
\end{align}
for  $i,j=1,2,\dots,q$. The prior distribution for $\bs{\beta}$ can then be written as
\begin{align}
    \pi(\bs{\beta})\propto\exp\left\{-\frac{1}{2}\norm[\GGamma_\beta^{-1}]{\bs{\beta}-\bs{\beta}_\ast}^2\right\}.
\end{align}

In the current paper, we take the boundary shape and Robin coefficient to be (statistically) mutually independent a priori, allowing us to decompose the prior as $\pi_{\rm prior}(\bs{\alpha},\bs{\beta})=\pi_{\alpha}(\bs{\alpha})\pi_{\beta}(\bs{\beta})$. The joint prior density is thus given by
\begin{align}
    \pi(\bs{\alpha},\bs{\beta})= \pi(\bs{\alpha})\pi(\bs{\beta})\propto\exp\left\{-\frac{1}{2}\norm[\GGamma_\alpha^{-1}]{\bs{\alpha}-\bs{\alpha}_\ast}^2-\frac{1}{2}\norm[\GGamma_\beta^{-1}]{\bs{\beta}-\bs{\beta}_\ast}^2\right\}.
\end{align}

\subsection{The Likelihood}
As is fairly standard, see e.g. \cite{NicholsonPetraKaipio18,Rundell08,CakoniKressSchuft10a}, we assume the data, $\data\in\mathbb{R}^m$, is corrupted by additive noise $\bs{e}$ and is related to the parameters through 
\begin{align}
\data=\fwd(\bs{\alpha},\bs{\beta})+\bs{e},
\end{align}
where $\fwd:\mathbb{R}^{2p+1}\times\mathbb{R}^{q}\rightarrow\mathbb{R}^m$ is termed the parameters-to-observable mapping\footnote{Evaluated at the functions corresponding to the Fourier coefficients $\bs{\alpha}$ and the finite element coefficient vector $\bs{\beta}$.}, and $\bs{e}\in\mathbb{R}^m$ denotes the additive noise. As discussed in Section \ref{sec: fwdprob}, inline with the literature, we take the data to be comprised of point-wise noisy measurements of $u$, satisfying (\ref{eq: fwd1}), at points along the accessible part of the boundary, $\Gamma_{\rm A}$. As such, we can rewrite the parameters-to-observable mapping as $\fwd(\bs{\alpha},\bs{\beta})=\bs{\mathcal{B}}\bs{u}\in\mathbb{R}^m$, where $\bs{u}={\rm vec}(u_1,u_2,\dots,u_l)$, and  $\bs{\mathcal{B}}$ denotes the linear observation (interpolation) operator which gives the values of each $u_i$ at the measurement locations. In practice, each of the $u_i$ are discretised and thus with a slight abuse of notation, we take the parameter-to-observable mapping to be of $\bs{\mathcal{B}}:\mathbb{R}^{r}\rightarrow\mathbb{R}^m$.

Inline with most literature, the additive noise is assumed to be normally distributed, i.e., $\bs{e}\sim\mathcal{N}(\bs{e}_\ast,\GGamma_e)$, and mutually independent of the parameters. These assumptions lead to the likelihood taking the form 
\begin{align}
\pi(\data\vert \bs{\alpha},\bs{\beta})\propto \exp\left\{-\frac{1}{2}\norm[\GGamma_e^{-1}]{\data-\bs{\mathcal{B}}\bs{u}}^2\right\},
\end{align}
see for example \cite{KaipioSomersalo05,CalvettiSomersalo07,Stuart10}.

\subsection{The Posterior density}
By employing the above parameterisations for the inaccessible part of the boundary, $\Gamma_{\rm I}$, and the Robin coefficient, $\beta$, the solution to the inference problem is the joint posterior density of $\bs{\alpha}$ and $\bs{\beta}$. This can be expressed using the finite dimensional version of Bayes' formula, as
\begin{align}
    \pi_{\rm post}(\bs{\alpha},\bs{\beta}\vert \data)&\propto\pi_{\rm like}(\data\vert \bs{\alpha},\bs{\beta})\pi_{\rm prior}(\bs{\alpha},\bs{\beta})\\
    &\propto \exp\left\{-\frac{1}{2}\left(\norm[\GGamma_e^{-1}]{\data-\bs{\mathcal{B}}\bs{u}}^2+\norm[\GGamma_\alpha^{-1}]{\bs{\alpha}-\bs{\alpha}_\ast}^2+\norm[\GGamma_\beta^{-1}]{\bs{\beta}-\bs{\beta}_\ast}^2\right)\right\}.
\end{align}

In general, if the parameter-to-observable mapping, $\fwd$, is nonlinear, the resulting posterior distribution is non-Gaussian. To explore the posterior distribution it is then standard to use sampling-based methods such as Markov chain Monte Carlo (MCMC) methods. However, it is well understood that these methods can be computationally expensive, especially when the parameter dimension is large, and/or if the forward problem itself is computationally expensive to run. To alleviate this problem (to some degree) we initially compute a Gaussian approximation to the posterior, the details of which are outlined in the following section. We define the posterior potential (of the posterior density),
\begin{align}
    \JJ:=\frac{1}{2}\left(\norm[\GGamma_e^{-1}]{\data-\bs{\mathcal{B}}\bs{u}}^2+\norm[\GGamma_\alpha^{-1}]{\bs{\alpha}-\bs{\alpha}_\ast}^2+\norm[\GGamma_\beta^{-1}]{\bs{\beta}-\bs{\beta}_\ast}^2\right),
\end{align}
so that the posterior can be written as $\pi_{\rm post}(\bs{\alpha},\bs{\beta}\vert \data)\propto\exp\left\{-\JJ\right\}$.

\subsubsection{The Laplace Approximation to the Posterior}
In the context of Bayesian inference for problems governed by partial differential equations (PDE), the most common Gaussian approximation to the posterior is the so-called Laplace approximation, see e.g.~\cite{KovalAlexanderianStadler20,NicholsonPetraKaipio18,Bui-ThanhGhattasMartinEtAl13}. Assuming the parameter-to-observable mapping is Fr\'{e}chet differentiable with respect to $\bs{\alpha}$ and $\beta$, the Laplace approximation is
$\hat{\pi}_{\rm post}(\bs{\alpha},\bs{\beta}\vert\data)=\mathcal{N}((\bs{\alpha}_{\rm MAP},\bs{\beta}_{\rm MAP}), \hat{\GGamma}_{\rm post})$, where $(\bs{\alpha}_{\rm MAP},\bs{\beta}_{\rm MAP})$ denotes the maximum a posteriori (MAP) estimate of $(\bs{\alpha},\bs{\beta})$, 
\begin{align}
    (\bs{\alpha}_{\rm MAP},\bs{\beta}_{\rm MAP})&:=\arg\max_{\bs{\alpha},\bs{\beta}}\pi_{\rm post}(\bs{\alpha},\bs{\beta})\\
    &\propto\arg\min_{\bs{\alpha},\bs{\beta}}\JJ,
\end{align}
and the approximate posterior covariance matrix is given by
\begin{align} \label{eq:Laplace-approximation}
    \hat{\GGamma}_{\rm post}=\left(\mathsf{G}(\bs{\alpha}_{\rm MAP},\bs{\beta}_{\rm MAP})^\ast\GGamma_e^{-1}\mathsf{G}(\bs{\alpha}_{\rm MAP},\bs{\beta}_{\rm MAP})+\GGamma_{\un}^{-1}\right)^{-1},
\end{align}
where $\GGamma_{\un}={\rm diag}(\GGamma_{\alpha},\GGamma_{\beta})$, and  $\mathsf{G}(\bs{\alpha}_{\rm MAP},\bs{\beta}_{\rm MAP})$ denotes the generalised derivative of $\mathcal{G}$ with respect to $\alpha$ and $\beta$ evaluated at the MAP estimate, i.e.,
\begin{align}\label{eq: Jaco}
    \mathsf{G}(\bs{\alpha}_{\rm MAP},\bs{\beta}_{\rm MAP})=[\mathsf{D}_{\alpha}\fwd[\bs{\alpha}_{\rm MAP},\bs{\beta}_{\rm MAP}] \quad \mathsf{D}_{\beta}\fwd[\bs{\alpha}_{\rm MAP},\bs{\beta}_{\rm MAP}]]
\end{align}

It's worth noting that the approximate posterior covariance matrix coincides with the inverse of the Gauss-Newton approximation to the Hessian of the negative log posterior (referred to simply as the Hessian in what follows), denoted \bs{H}, i.e.,
\begin{align}
   \hat{\GGamma}_{\rm post}=\bs{H}^{-1}. 
\end{align}

\subsubsection{Full Characterisation of the Posterior}
In some cases, the MAP estimate and a localised Gaussian approximation, such as the Laplace approximation, may be sufficient for the quantification of uncertainty in the parameters.
However, without fully characterising the posterior it is essentially impossible to determine the feasibility, and thus relevance, of the local Gaussian approximation.
As stated previously, fully characterising the posterior for nontrivial problems requires sampling based methods such as MCMC.

It is well understood that MCMC samplers can become practically infeasible for high dimensional problems.
However, this {\em curse of dimensionality} can often be offset to some extent by implementing sampling schemes which exploit the geometry of the posterior, see for example \cite{GirolamiCalderhead11,BuiGirolami14}.
One of the most popular MCMC sampling schemes which takes into account the geometry of the posterior is the Metropolis-adjusted Langevin algorithm (MALA) \cite{RobertsTweedie96,roberts1998optimal,ApteHairerStuartEtAl07}, which can be motivated using Langevin dynamics.
The basic idea in MALA is to guide the proposal towards an area of higher probability, thus speeding up the convergence.
Use of MALA requires the gradient of the negative log posterior with respect to the unknown parameters.
By employing the so-called adjoint approach \cite{BorziSchulz11,HinzePinnauUlbrichEtAl08,Troltzsch10}, the gradient can be efficiently calculated by solving one additional (adjoint-) forward solve, see Section \ref{sec: adjs} below.

We now briefly recall the key steps to MALA, for more details see for example~\cite{RobertsTweedie96,roberts1998optimal}.
Langevin diffusion is defined by the following stochastic differential equation
\begin{align}\label{eq: SDE}
    d{\bs{m}}_t=-\bs{A}\nabla_{\bs{m}} \JJ(\bs{m}) dt +\sqrt{2}\bs{A}^{1/2}\,d\bs{B}_t,
\end{align}

\noindent where we have denoted all the unknowns as $\bs{m}=(\bs{\alpha},\bs{\beta})$, $\bs{A}$ is a positive definite (preconditioning) matrix, and $\bs{B}_t$ is standard Brownian motion. 

A key feature of the above process is that in the continuous-time case its stationary and limiting distribution is $\pi_{\rm post}$.
However, discrete approximations to \eqref{eq: SDE} can have vastly different asymptotic behaviours from the diffusion process they try to approximate \cite{RobertsTweedie96}.
As such, it is necessary to introduce a Metropolis-Hastings accept/reject step that ensures the convergence to $\pi_{\rm post}$.
In MALA, the proposal for the next step, $\un^\ast_{k+1}$,  is generated as
\begin{align}
   {\bs{m}}^\ast_{k+1}=\un_{k}-\tau\bs{A}\nabla_{\bs{m}} \JJ(\bs{m}_k) + \sqrt{2\tau}\bs{A}^{1/2}\,\bs{\xi},
\end{align}

\noindent where now the preconditioning matrix $\bs{A}$ represents the covariance of the proposal density $q(\bs{m}_{k+1}^{\ast} \vert \bs{m}_k)$, $\tau>0$ can be seen as a step size parameter that globally scales the proposal density, and $\bs{\xi}$ is a standard normal ($q+2p+1$-dimensional) probability density.
The proposal is then accepted or rejected according to the Metropolis-Hastings algorithm, where
\begin{align} \label{eq:accratio}
   \alpha = \min\left\{1, \frac{\pi(\un^\ast_{k+1}) q(\un_k\vert\un^\ast_{k+1})}{\pi(\un_k) q(\un^\ast_{k+1}\vert\un_k)}\right\},
\end{align}

\noindent is the probability of accepting the proposal.
In practice, we work with the logarithm of $\alpha$ to avoid numerical underflow.

In MALA, the drift towards the posterior gradient makes the proposal density $q(\cdot\vert\cdot)$ nonsymmetric, and therefore the ratio $q(\un_k\vert\un^\ast_{k+1})/q(\un^\ast_{k+1}\vert\un_k)$ does not cancel out.
The log-ratio of the proposal densities is given by
\begin{align}
	\log\left(\frac{q(\un_k\vert\un^\ast_{k+1})}{q(\un^\ast_{k+1}\vert\un_k)}\right) = -\frac{1}{4\tau}\Big(\norm{\bs{L}(-\delta\bs{m} + \tau\bs{A}\nabla_{\bs{m}} \JJ(\bs{m}_{k+1}^\ast))}^2
	- \norm{\bs{L}(\delta\bs{m} + \tau\bs{A}\nabla_{\bs{m}} \JJ(\bs{m}_{k}))}^2\Big),
\end{align}
where $\delta\bs{m} := \bs{m}_{k+1}^\ast - \bs{m}_{k}$ is the taken step, and $\bs{L}^T\bs{L} = \bs{A}^{-1}$.

Let us now discuss the selection of the proposal covariance $\bs{A}$ and the proposal scaling parameter $\tau$.
To achieve optimal sampling efficiency, we should have $\bs{A} = \bs{\Sigma}_\pi$, where $\bs{\Sigma}_\pi$ denotes the covariance of the true posterior \cite{atchade2006adaptive}.
In addition, it has been shown that under various assumptions the asymptotically optimal acceptance rate for MALA is 0.574 \cite{roberts1998optimal}.
The acceptance rate can be adjusted using the scaling parameter~$\tau$.
Since $\bs{\Sigma}_\pi$ and the optimal value of $\tau$ are obviously unknown prior to doing the inversion, we adapt them continuously during the sampling, see \cite{HaarioSaksmanTamminen01, atchade2006adaptive, AndrieuThoms08}.
A reasonable starting location for MCMC is the MAP estimate given by the Gauss-Newton iterations, and as the initial proposal covariance before the start of the adaptation we use the Laplace approximation \eqref{eq:Laplace-approximation}.

Before any conclusions can be drawn from a MCMC run, we need to be reasonably sure that the chain has converged to its limiting distribution.
Naturally, the convergence can never be guaranteed since any method of diagnosing convergence can by definition only diagnose the non-convergence of a chain.
In practice, however, we can arrive at a reasonable certainty by testing the convergence with multiple complementary methods.
These methods include visual inspection of the chains and measuring the sampling quality by autocorrelation, comparing within-chain variances between multiple runs (the Gelman-Rubin diagnostic) \cite{gelman1992inference,brooks1998general}, and estimating sampling uncertainty by computing the Monte Carlo standard error (MCSE) \cite{flegal2008markov, vats2019multivariate}.

In this paper, we implement a stopping rule based on the estimated MCSE and its relation to the posterior variance of the parameters, as in \cite{NiskanenDazelGrobyEtAl19}.
First, we run the sampler for enough steps so that the MCSE can be reliably calculated using the non-overlapping batch means method, see \cite{flegal2008markov}.
Then, at regular intervals we check, for each unknown individually, if the MCSE estimate at 98~\% confidence level is smaller than 10~\% of the posterior standard deviation.
When this is true for each unknown, we terminate the run.

\subsection{The Adjoint Approach and Further Derivatives}\label{sec: adjs}
Calculation of the MAP estimate, $(\bs{\alpha}_{\rm MAP},\bs{\beta}_{\rm MAP})$, and construction of the Laplace approximation of the posterior covariance matrix, $\hat{\GGamma}_{\rm post}$, requires the generalised derivative $\mathsf{G}$ (see \ref{eq: Jaco}), details of which were provided in Section \ref{sec: Der1}. On the other hand, the use of MALA necessitates computation of only the gradient of the negative log posterior with respect to both the inaccessible part of the boundary and the Robin coefficient. Although the gradient is essentially trivial to calculate with the Jacobian in hand, computing the Jacobian can be costly, requiring as many forward solves as the number of measurements. To avoid this potential bottleneck, we employ the so-called adjoint approach~\cite{BorziSchulz11,HinzePinnauUlbrichEtAl08,Troltzsch10}, which computes the gradient at a cost of one additional (adjoint) forward solve. In what follows, we will set $g=g_k$ and $u=u_k$ for some fixed $k\in\left\{1,2,\dots,\ell\right\}$. We begin by introducing  a Lagrangian, $\mathcal{L}:H^1_D(\Omega)\times H^1_D(\Omega)\times \mathbb{R}^{2p+1}\times H^1([0,L])\rightarrow \mathbb{R}$, as
\begin{align}
    \mathcal{L}(u,v,\bs{\alpha},\beta)&=\JJ+\int_{\Omega(\bs{\alpha})}\sigma\nabla u\cdot\nabla v\;d\bs{x}+\int_{\Gamma_{\rm I}(\bs{\alpha})}\exp{(\beta)} u v\;d\bs{s}_{\rm I}-\int_{\Gamma_{\rm A}}g v\;d\bs{s}_{\rm A}\\
    &=\JJ+\int_{\oO}\os_{\bs{\alpha}}\nabla \ou\cdot\nabla \tilde{v}\;d\tilde{\bs{x}}+\int_{\oG_{\rm I}}\exp{(\oB_{\bs{\alpha}})} \ou \tilde{v}\;d\tilde{\bs{s}}_{\rm I}-\int_{\Gamma_{\rm A}}g \tilde{v}\;d\bs{s}_{\rm A},
\end{align}
where we have used the invariance property of the Poisson equation, and where the subscript $\bs{\alpha}$ is used to denote the dependence on $\bs{\alpha}$. 

Determining the gradient of $\JJ$ (with respect to the unknowns) is achieved by requiring that variations of the Lagrangian with respect to the
forward potential, $u$, and the so-called adjoint potential, $v$, are 0, see e.g. \cite{Troltzsch10}. The variation of the Lagrangian with respect to the adjoint potential gives \begin{align}
   \mathsf{D}_{\tilde{w}}\mathcal{L}[u,v,\bs{\alpha},\beta]=\JJ+\int_{\oO}\os_{\bs{\alpha}}\nabla \tilde{u}\cdot\nabla \tilde{w}\;d\tilde{\bs{x}}+\int_{\oG_{\rm I}}\exp{(\oB_{\bs{\alpha}})} \tilde{u} \tilde{w}\;d\tilde{\bs{s}}_{\rm I}-\int_{\Gamma_{\rm A}}g \tilde{w}\;d\bs{s}_{\rm A}=0,
\end{align}
for $\tilde{w}\in H^1_D(\oO)$, i.e., $u$ must satisfy the forward problem. The variation with respect to the forward potential give the adjoint equation,
\begin{align}
   \mathsf{D}_{\tilde{z}}\mathcal{L}[u,v,\bs{\alpha},\beta]=-\bs{\mathcal{B}}^\ast\GGamma_e^{-1}(\data-\bs{\mathcal{B}}\tilde{z})+\int_{\oO}\os_{\bs{\alpha}}\nabla \tilde{z}\cdot\nabla \tilde{v}\;d\tilde{\bs{x}}+\int_{\oG_{\rm I}}\exp{(\oB_{\bs{\alpha}})} \tilde{z} \tilde{v}\;d\tilde{\bs{s}_{\rm I}}=0
\end{align}
for $\tilde{z}\in H^1_D(\oO)$, which the adjoint variable $v$ must satisfy (this is the one additional forward solve required). Finally, the gradients with respect to the Fourier coefficients and the Robin coefficient are then given (in strong form) as 
\begin{align}
  \nabla_{\bs{\alpha}}\mathcal{L}(u,v,\bs{\alpha},\beta)=\GGamma_{\alpha}^{-1}(\bs{\alpha}-\bs{\alpha}_\ast)+\int_{\oO}\bs{\kappa}_{\bs{\alpha}}\nabla \tilde{u}\cdot\nabla \tilde{v}\;d\tilde{\bs{x}}+\int_{\oG_{\rm I}}{\bs b}_{\bs{\alpha}}\tilde{u} \tilde{v}\;d\tilde{\bs{s}}_{\rm I},
\end{align}
and
\begin{align}
   \nabla_\beta\mathcal{L}(u,v,\bs{\alpha},\beta)=\mathcal{A}_{\beta}(\beta_{\bs{\alpha}}-\beta_\ast)+\exp{(\beta_{\bs{\alpha}})}u v,
\end{align}
respectively, where 
\begin{align}
\left[\bs{\kappa}_{\bs{\alpha}}\right]_i(\tilde{\bs{x}})=\sigma(\bs{x})\frac{\partial}{\partial\alpha_i}\begin{pmatrix} 
\begin{array}{c|c}
f(x_1) & -\frac{x_2}{f(x_1)}\frac{df(x_1)}{dx_1} \\ 
\hline
-\frac{x_2}{f(x_1)}\frac{df(x_1)}{dx_1} & \frac{1}{f(x_1)}+\frac{x_2^2}{f^3(x_1)}\left(\frac{df(x_1)}{dx_1}\right)^2
\end{array}
\end{pmatrix},\quad
\left[{\bs b}_{\bs{\alpha}}\right]_i(\tilde{s})=\exp\left(\beta(s)\right)\frac{\partial}{\partial\alpha_i}\sqrt{1+\left(\frac{df}{ds}\right)^2H^2},
\end{align}
where ${\bs{x}=\psi^{-1}(\tilde{\bs{x}})}$ and ${s=\varphi^{-1}(\tilde{s})}$, and $i=0,1,\dots, 2p$.

Finally, the gradient of the cost function evaluated at $(u,v,\bs{\alpha},\beta)$ is then given by 
\begin{align}
\nabla\JJ=\begin{pmatrix}\nabla_{\bs{\alpha}}\mathcal{L}(u,v,\bs{\alpha},\beta)\\
\nabla_{\beta}\mathcal{L}(u,v,\bs{\alpha},\beta)
\end{pmatrix}.
\end{align}

\section{Numerical Examples}\label{sec: numex}
In this section, we outline three numerical examples to assess the applicability, performance, and robustness of the
proposed approach.
In each example, we take the uncorroded (reference) domain to be a slab with  with length $L=1$, and height $H=0.05$.

To reduce so-called {\em inverse crimes} \cite{KaipioSomersalo07}, the data-generating meshes have a substantially finer discretisation than the meshes used at the inference stage. Furthermore, in each example, neither the true Robin coefficient $\beta$, nor the true boundary shape, $\Gamma_{\rm I}$, are generated form their respective priors, which we describe below.

The synthetic data is generated solving (\ref{eq: fwd1}) using the finite element method using $\ell = 8$ different Neumann boundary conditions (see (\ref{eq: Neumann})). Measurement data consists of 32 equally spaced noisy point-wise measurements of $u_\ell$ along $\Gamma_{\rm A}$, and the as shown in Figure~\ref{fig:fwdmeshes} using black dots.
In each example the noise is taken to have mean $\bs{e}_\ast=\bs{0}$ and covariance matrix $\GGamma_e=\delta_e^2\bs{I}$, where $\delta_e=\left(\max(\mathcal{B}u)-\min(\mathcal{B}u)\right) \times 1/100$, that is, the noise level is $1\%$ of the range of the noiseless measurements.
\begin{itemize}
\item
{\bf Example 1:} The first example we consider exhibits a smooth, shallow dip in the inaccessible part of the boundary (see Figure~\ref{fig:fwdmeshes}), and a smooth true Robin coefficient.
For brevity, the true Robin coefficients for all examples are shown with the inversion results in Figures~\ref{fig:results_case1}-\ref{fig:results_case3}.
\item
{\bf Example 2:} The second example has a wider corroded area on the left side of the slab, and an increase in the thickness at the right, which could represent an area with deposition or buildup. Furthermore, changes in the boundary are steeper relative to Example 1. In addition, a small amount of white noise is added to the boundary shape and the Robin coefficient to model, perhaps, a more realistic situation where the boundary and Robin coefficient are {\em rough}.
\item
{\bf Example 3:} Finally, in the third example we place three deep and steep cavities across the inaccessible part of the slab.
This example is the most challenging as the true boundary shape and Robin coefficient have 0 prior probability. However, this example does provide a good test of the problem's robustness to the choice of prior.
\end{itemize}

Before discussing the results of each example we outline further computational details.
\subsection{Computational Details}

\begin{figure}[t]
	\includegraphics[width=\linewidth]{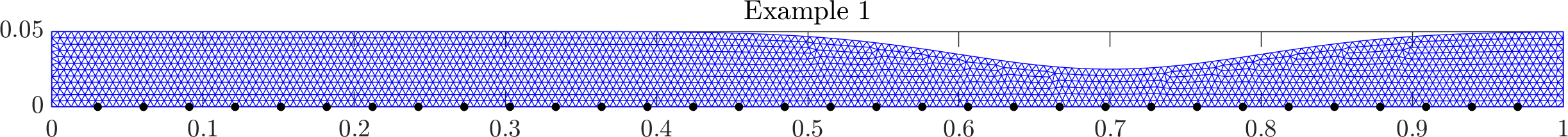}
	\vspace{0.1mm}
	
	\includegraphics[width=\linewidth]{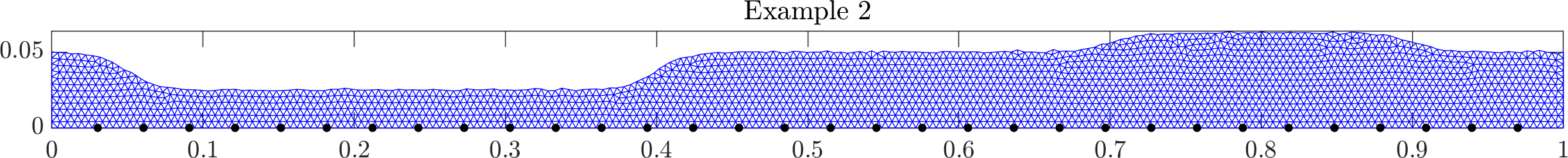}
	\vspace{0.1mm}
	
	\includegraphics[width=\linewidth]{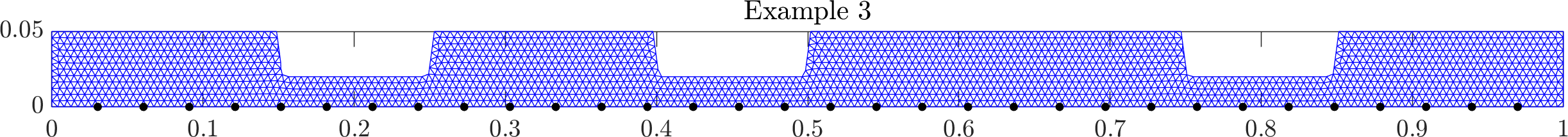}
	\caption{Data-generating meshes of the 3 examples considered. The black dots denote the measurement locations.}
	\label{fig:fwdmeshes}
\end{figure}

For each example considered, we take the diffeomorphism (from the true domain $\Omega$ to $\oO$) to be of the form
\begin{align}
    \psi:(x_1,x_2)\mapsto\left(x_1,x_2(f(x_1))^{-1}\right), \quad (x_1,x_2)\in\Omega,
\end{align}
with
\begin{align}
   f(x_1) = 1 + \alpha_0 + \sum_{n=1}^p\alpha_{2n} \cos(2\pi nx_1)+\alpha_{2n-1} \sin(2\pi nx_1).
\end{align}
The choice of $p$ can be viewed as making a trade-off between computational speed and the maximum possible accuracy of the solution. That is, with a smaller $p$ we have fewer parameters to infer but smaller values of $p$ cannot represent as complex shapes. In the current paper we take $p = 7$, meaning we have 15 Fourier coefficients to estimate. 

\begin{figure}[t]
	\includegraphics[width=\linewidth]{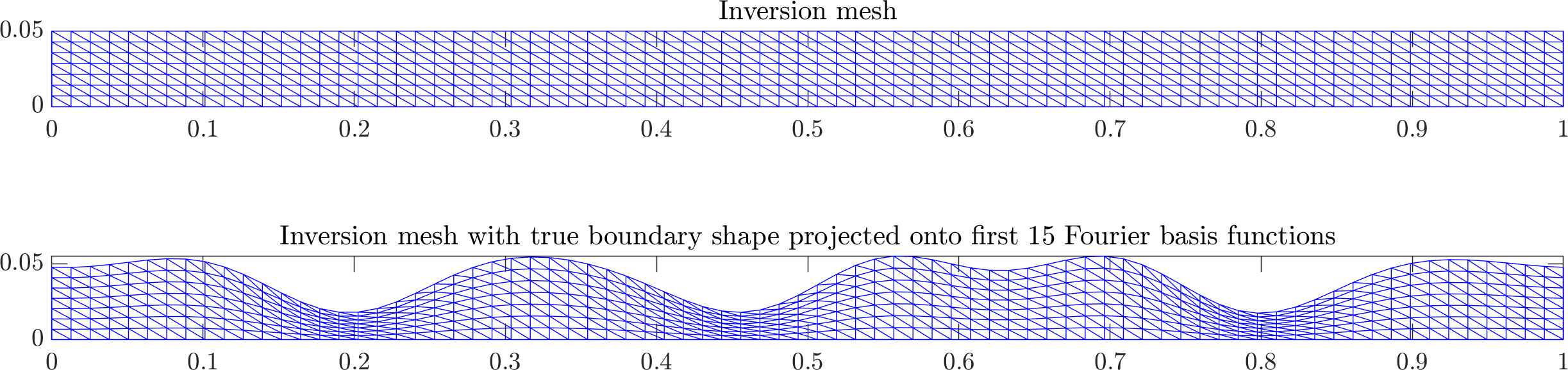}
	\caption{Top: The reference inversion mesh used for all examples. Bottom: True domain discretised using the inversion mesh with the true boundary shape projected onto the first 15 Fourier basis functions.}
	\label{fig:inversion_mesh_example}
\end{figure}

For Examples 2 and 3, representing the true boundary shapes accurately requires substantially more basis functions than the 15 used.
Thus, to better assess the estimates, we can compare them to the true boundary shape projected onto the truncated Fourier basis.
This is illustrated in Figure~\ref{fig:inversion_mesh_example}, where we show the reference inversion mesh (used in all examples) as well as the inversion mesh corresponding to the true domain for Example 3 using the first 15 Fourier basis functions. It's clear that steeper changes in the boundary shape cannot be represented by the truncated Fourier series.

\begin{figure}[t]
	\includegraphics[width=\linewidth]{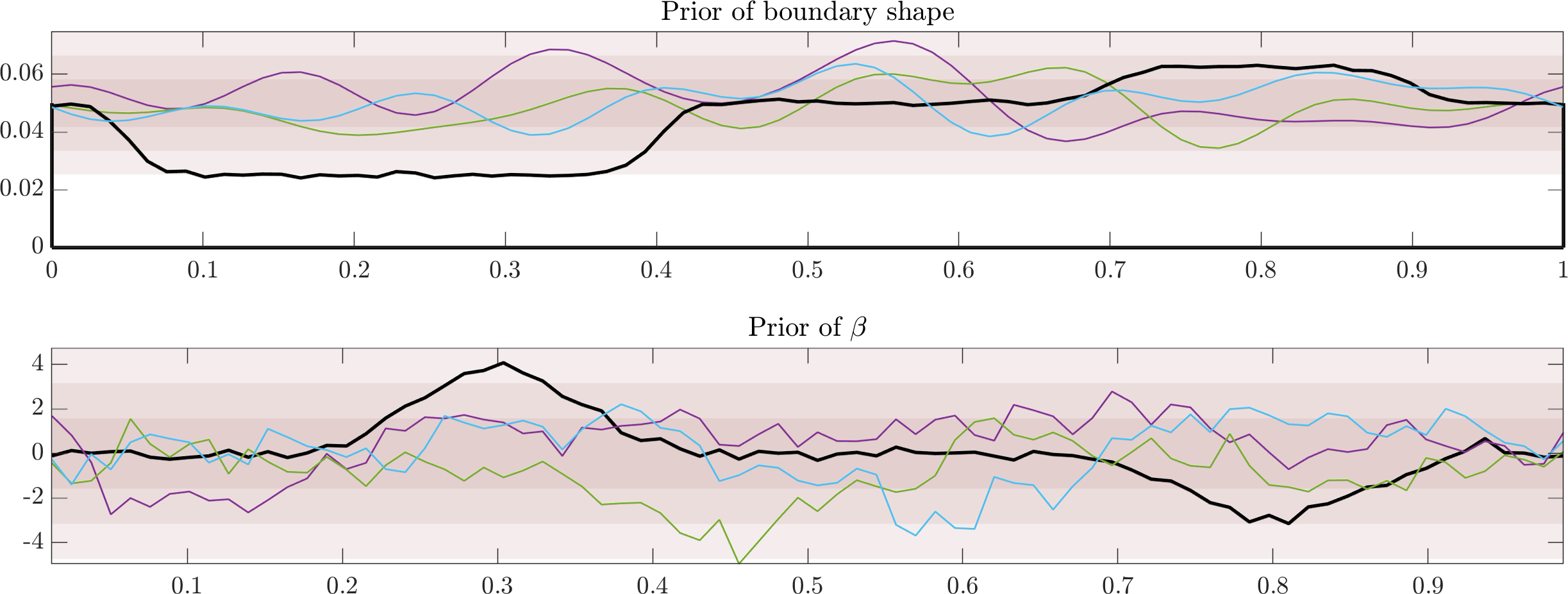}
	\caption{Draws from the prior over the true values of Example 2. The shaded areas indicate the one to three standard deviations of the prior.}
	\label{fig:prior_draws}
\end{figure}

We use the same prior distribution for all examples.
Specifically, for the boundary shape we take $\bs{\alpha}_\ast = \bs{0}$ and
\begin{equation}
	\GGamma_\alpha =\mathrm{diag}(\sigma_0^2,\sigma_1^2,\sigma_1^2,\dots,\sigma_p^2,\sigma_p^2)=\sigma_\alpha^2 \mathrm{diag}(1, 2^{s_\alpha}, 2^{s_\alpha}, \dots, (p+1)^{s_\alpha}, (p+1)^{s_\alpha}),
\end{equation}
where $\sigma_\alpha^2$ controls the overall variance of the parameters, and $s_\alpha<0$ controls how quickly the variance diminishes towards the higher frequency components.
The faster the Fourier coefficients converge towards zero, the smoother the resulting function is \cite{Grafakos08}, and as such, $s_\alpha$ can be interpreted as a smoothness parameter. We choose $\sigma_\alpha^2 = 0.01$ and $s_\alpha = -1$ for all considered examples. For the Robin coefficient we set $\bs{\beta}_\ast = \bs{0}$ and $\GGamma_\beta$ as in \eqref{eq:beta_covariance}, where we set $\delta_\beta^2 =50$ and $l=10$.
Draws from both priors are shown in Figure~\ref{fig:prior_draws}.

The meshes used for the generation of the synthetic data are different in each case (see Figure~\ref{fig:fwdmeshes}), and each mesh has between 2,000 and 3,000 nodes. On the other hand, for all examples, in the data-generation mesh the Robin coefficient is represented using 230 finite element basis functions. The discretisation of the mesh used for inference has a total of 640 nodes, and results in 78 degrees of freedom for the Robin coefficient, i.e., $\bs{\beta}\in\mathbb{R}^{78}$.
Thus, the total number of the parameters to be inferred is $78 + 15 = 93$.

Finally, we remark that, inline with e.g.~\cite{EisenstatWalker96}, the Gauss-Newton method used to compute the MAP estimate is terminated when the norm of the gradient has decreased by a factor of $10^5$, while {\em line search} is carried out using the Armijo condition with $c_{1}=10^{-4}$ and $c_{2}=0.9$, see ~\cite{NocedalWright06}. All computations were carried out using MATLAB 2018b on a laptop with an Intel i7-8850H CPU and 16 GB RAM.


\subsection{Results} \label{sec:results}
Here we discuss and compare the Laplace approximation of the posterior and the accurate posterior found using MCMC for each of the examples outlined above. 

\begin{figure}[h]
	\includegraphics[width=\linewidth]{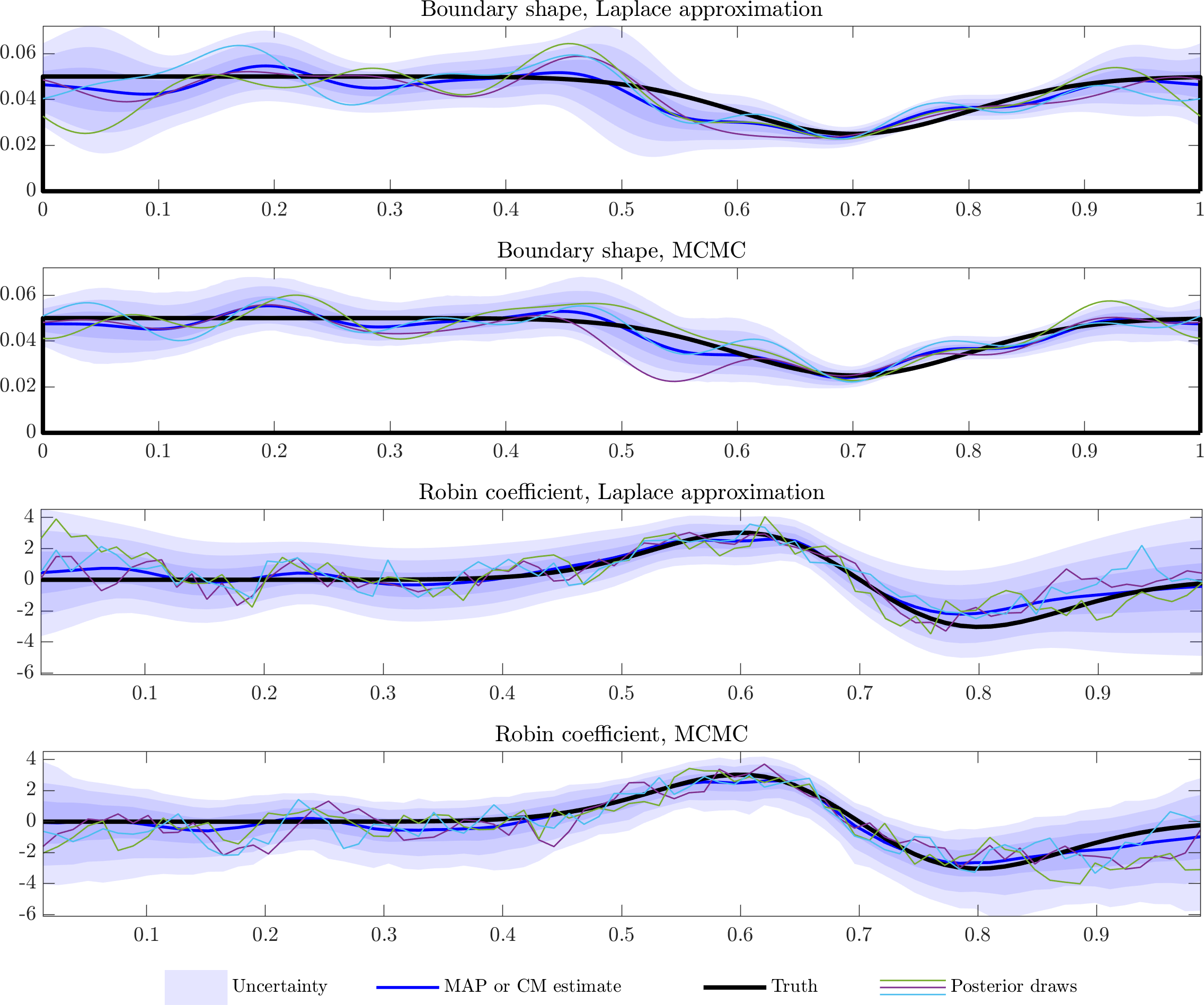}
	\caption{Results for Example 1: MAP and CM estimates of the boundary shape and Robin coefficient, some posterior draws, and the true boundary shape and Robin coefficients. The shaded areas denote the posterior uncertainty estimates as one, two, and three sigma. These are computed as MAP-estimate $\pm \:\: n\sigma, n=1,2,3$ for the Laplace approximation, and as $68-95-99.7$ credibility regions for the MCMC results.}
	\label{fig:results_case1}
\end{figure}

{\bf Example 1.} For the first example, the MAP estimate was found after 24 Gauss-Newton iterations, while the MCMC convergence criterion was reached in approximately 70,000 steps after burn-in (10,000 steps in all examples), which took approximately two hours. The MAP estimate, Laplace approximation of the (marginal) posterior, and the conditional mean and the (marginal) posterior found using MCMC are all shown in Figure~\ref{fig:results_case1}.

Both the Laplace approximation of the posterior and accurate posterior support the true boundary shape and Robin coefficient well. Comparing the Laplace approximation and accurate posterior for the boundary shape, it is clear that the accurate posterior has significantly narrower confidence intervals, particularly towards the lateral boundaries. Moreover, for the boundary shape, the CM and MAP estimates are in fairly good agreement, while the accurate posterior seems to be symmetric, indicating that the posterior is well approximated by a normal distribution. On the other hand, when comparing the Laplace approximation and accurate posterior for the Robin coefficient, we see that the accurate posterior is negatively skewed (evident from the confidence intervals and the fact that the CM estimate is generally smaller than the MAP estimate). 

\begin{figure}[h]
	\includegraphics[width=\linewidth]{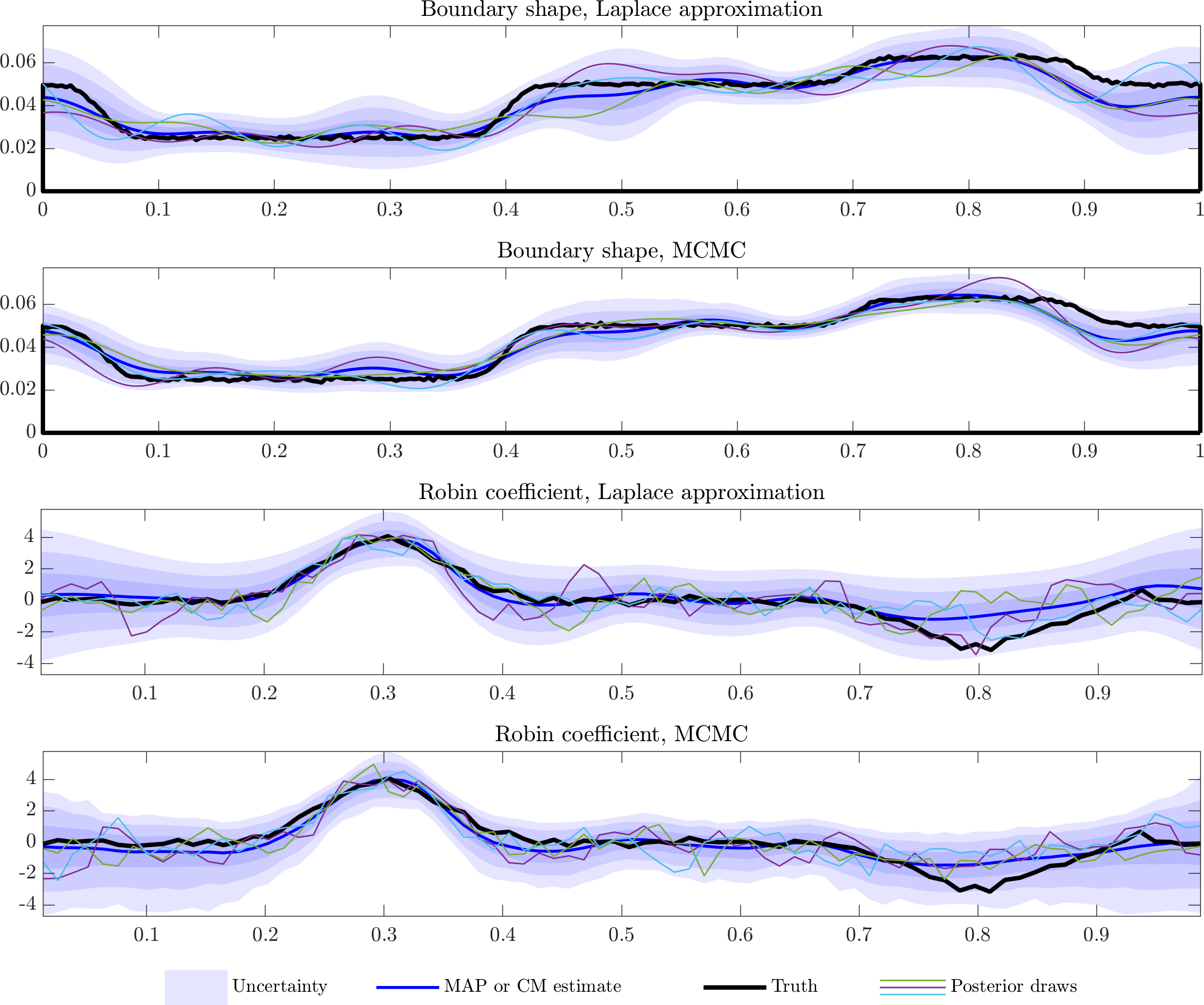}
	\caption{Results for Example 2.}
	\label{fig:results_case2}
\end{figure}

{\bf Example 2.} For the second example, calculating the MAP estimate took 30 Gauss-Newton iterations, while the MCMC convergence criterion was reached in approximately 140,000  steps after burn-in, which took approximately three and a half hours. In Figure~\ref{fig:results_case2} we show the MAP estimate, Laplace approximation of the posterior and the posterior and conditional mean found using MCMC. 

As in Example 1, the Laplace approximation of the posterior and accurate posterior support the true boundary shape and Robin coefficient well. It also clear that for both the boundary shape and Robin coefficient, estimates are fairly insensitive to the the small, highly oscillatory, variations in corresponding truths. Similarly to Example 1, the Laplace approximation has slightly wider confidence intervals compared to the accurate posterior, particularly towards the lateral boundaries, and the CM and MAP estimates are in good agreement for the boundary shape.
Again, the accurate posterior for the Robin coefficient is negatively skewed. 

\begin{figure}[h]
	\includegraphics[width=\linewidth]{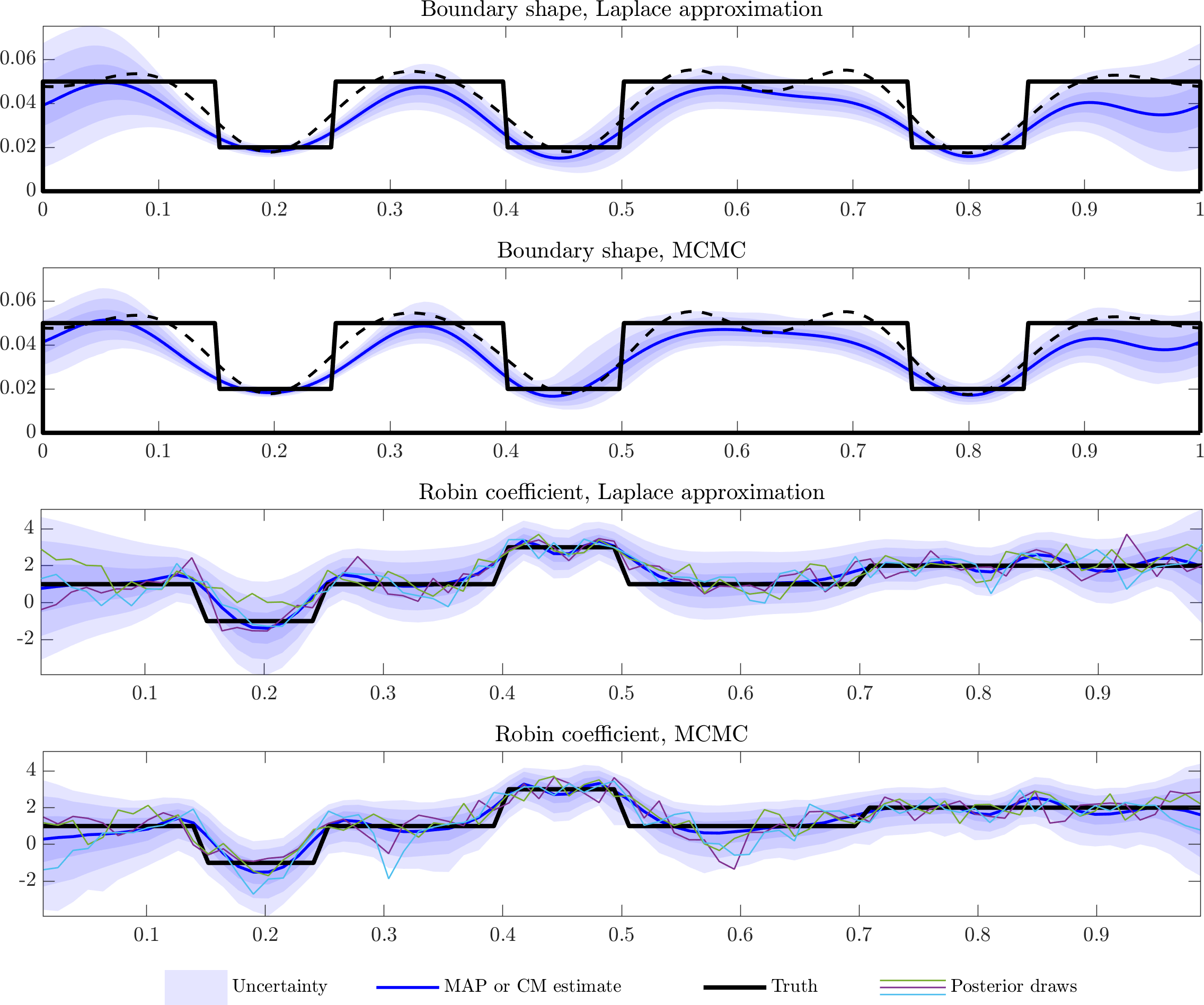}
	\caption{Results for Example 3. The black dashed line shows the projection of the true boundary shape onto the first 15 Fourier basis coefficients.}
	\label{fig:results_case3}
\end{figure}

{\bf Example 3.}  For the final example 31 Gauss-Newton iterations were needed to compute the MAP estimate, while the MCMC convergence criterion was reached in approximately 80,000  steps after burn-in, taking approximately two hours. Figure~\ref{fig:results_case3} shows the MAP estimate, Laplace approximation of the posterior and the posterior and CM found using MCMC. Since the true boundary shape lies well outside the prior distribution, we also show the projection of the true boundary shape onto the first 15 Fourier basis functions. This provides a reference for the best possible reconstruction.

As in the previous examples, the Laplace approximation to the posterior and the accurate posterior for the Robin coefficient support the true coefficient well. Conversely, though as expected, the sharp jumps in the boundary shape are poorly supported by both the approximate and accurate posterior densities as these are impossible to estimate using the  given prior on the boundary shape. However, we do see that the both posterior densities for the boundary shape support the projection of the true boundary shape onto the first 15 Fourier basis functions fairly well. It is also evident that the uncertainty in the boundary shape is significantly lower in regions where the domain is narrowest. This is inline with our intuition as the data (collected on the bottom of the slab) is much more informative on the shape of the boundary which is closer to where the data is measured. Finally, the accurate posterior does show some improvement over the approximate posterior for the boundary shape, particularly at the sharp jumps, as well as being slightly positively skewed. 

\section{Conclusion}
In this paper, we considered the simultaneous inference of both the shape of an inaccessible part of the domain boundary, as well as the Robin coefficient on the inaccessible part of the domain boundary. We considered three examples with different shapes and associated Robin boundary conditions for the inaccessible part of the domain. By exploiting an invariance property of the Poisson problem we were able to avoid having to re-mesh the changing domain shape, which in turn negates the need to recompute the discretised covariance operator of the Robin coefficient on the inaccessible part of the domain. Furthermore, use of the invariance property allowed for a straightforward and efficient means of calculating the required shape derivatives. 

In each case, we carried out both local (optimisation-based), approximate, uncertainty quantification in both parameters, as well as full exploration of the posterior using MCMC with MALA. To ensure efficiency of the MCMC algorithm, and possible application to larger-scale problems, we initially posed
the problem in infinite dimensions and then employed the
adjoint approach to compute the required derivatives.

Our results suggest that for the specific cases of the problem considered here, the Laplace approximation generally provides a fairly accurate representation of the true posterior, though it cannot model any skew, which was present in the accurate posterior densities found using MCMC. However, in assessing the applicability and performance of our proposed approach the current study was only carried out in two dimensions, with domain shape changes essentially restricted to the top surface only. Furthermore, the types of boundary shapes which could be estimated was limited to those which could be represented using only a small number of Fourier basis functions. Finally, though common in the literature, the measurement setup used is hard to carry out in reality. Natural next steps for future work then, would be to apply the same framework to three dimensional problems while considering a more realistic measurement set up as well as investigating other types of representations for the shape of the inaccessible part of the domain.

\section{Acknowledgements}
The authors wish to acknowledge Tom ter Elst and Gareth Gordon for their insights into the use of diffeomorphisms as well as Jari Kaipio for several valuable discussions.

\bibliographystyle{unsrt}
 \bibliography{DiffeoBib}

\begin{thebibliography}{10}

\bibitem{KaupSantosa95}
Peter~G Kaup and Fadil Santosa.
\newblock Nondestructive evaluation of corrosion damage using electrostatic
  measurements.
\newblock {\em Journal of Nondestructive Evaluation}, 14(3):127--136, 1995.

\bibitem{Inglese97}
Gabriele Inglese.
\newblock An inverse problem in corrosion detection.
\newblock {\em Inverse problems}, 13(4):977, 1997.

\bibitem{AlessandriniDelPieroRondi03}
Giovanni Alessandrini, L~Del~Piero, and Luca Rondi.
\newblock Stable determination of corrosion by a single electrostatic boundary
  measurement.
\newblock {\em Inverse problems}, 19(4):973, 2003.

\bibitem{DivoKassabKapatEtAl05}
Eduardo Divo, Alain~J Kassab, Jay~S Kapat, and Ming-King Chyu.
\newblock Retrieval of multidimensional heat transfer coefficient distributions
  using an inverse {BEM}-based regularized algorithm: numerical and
  experimental results.
\newblock {\em Engineering Analysis with Boundary Elements}, 29(2):150--160,
  2005.

\bibitem{Shantasiriwan99}
Somchart Chantasiriwan.
\newblock Inverse heat conduction problem of determining time-dependent heat
  transfer coefficient.
\newblock {\em International journal of heat and mass transfer},
  42(23):4275--4285, 1999.

\bibitem{ChaabaneJaoua99}
Slim Chaabane and Mohamed Jaoua.
\newblock Identification of {R}obin coefficients by the means of boundary
  measurements.
\newblock {\em Inverse problems}, 15(6):1425, 1999.

\bibitem{NicholsonPetraKaipio18}
Ruanui Nicholson, No{\'e}mi Petra, and Jari~P Kaipio.
\newblock Estimation of the {R}obin coefficient field in a {P}oisson problem
  with uncertain conductivity field.
\newblock {\em Inverse Problems}, 34(11):115005, 2018.

\bibitem{FasinoInglese99}
Dario Fasino and Gabriele Inglese.
\newblock An inverse {R}obin problem for {L}aplace's equation: {T}heoretical
  results and numerical methods.
\newblock {\em Inverse problems}, 15(1):41, 1999.

\bibitem{Jin07}
Bangti Jin.
\newblock Conjugate gradient method for the {R}obin inverse problem associated
  with the {L}aplace equation.
\newblock {\em International journal for numerical methods in engineering},
  71(4):433--453, 2007.

\bibitem{BaratchartBourgeoisLeblond16}
Laurent Baratchart, Laurent Bourgeois, and Juliette Leblond.
\newblock Uniqueness results for inverse {R}obin problems with bounded
  coefficient.
\newblock {\em Journal of Functional Analysis}, 270(7):2508--2542, 2016.

\bibitem{TossavainenVauhkonenHeikkinenEtAl04}
Olli-Pekka Tossavainen, M~Vauhkonen, LM~Heikkinen, and T~Savolainen.
\newblock Estimating shapes and free surfaces with electrical impedance
  tomography.
\newblock {\em Measurement Science and Technology}, 15(7):1402, 2004.

\bibitem{AlbuquerqueLaurainSturm20}
Yuri~Flores Albuquerque, Antoine Laurain, and Kevin Sturm.
\newblock A shape optimization approach for electrical impedance tomography
  with point measurements.
\newblock {\em Inverse Problems}, 2020.

\bibitem{BanksKojimaWinfree90}
H~Thomas Banks, Fumio Kojima, and WP~Winfree.
\newblock Boundary estimation problems arising in thermal tomography.
\newblock {\em Inverse problems}, 6(6):897, 1990.

\bibitem{FanSunYang08}
Chunli Fan, Fengrui Sun, and Li~Yang.
\newblock A new computational scheme on quantitative inner pipe boundary
  identification based on the estimation of effective thermal conductivity.
\newblock {\em Journal of Physics D: Applied Physics}, 41(20):205501, 2008.

\bibitem{HettlichRundell96}
Frank Hettlich and William Rundell.
\newblock Iterative methods for the reconstruction of an inverse potential
  problem.
\newblock {\em Inverse problems}, 12(3):251, 1996.

\bibitem{CakoniKress07}
Fioralba Cakoni and Rainer Kress.
\newblock Integral equations for inverse problems in corrosion detection from
  partial {C}auchy data.
\newblock {\em Inverse Problems \& Imaging}, 1(2):229, 2007.

\bibitem{Bacchelli08}
Valeria Bacchelli.
\newblock Uniqueness for the determination of unknown boundary and impedance
  with the homogeneous {R}obin condition.
\newblock {\em Inverse Problems}, 25(1):015004, 2008.

\bibitem{FasinoIngleseMariani07}
Dario Fasino, Gabriele Inglese, and Francesca Mariani.
\newblock Corrosion detection in conducting boundaries: {II}. {L}inearization,
  stability and discretization.
\newblock {\em Inverse Problems}, 23(3):1101, 2007.

\bibitem{IngleseMariani04}
G~Inglese and F~Mariani.
\newblock Corrosion detection in conducting boundaries.
\newblock {\em Inverse Problems}, 20(4):1207, 2004.

\bibitem{HarrachMeftahi19}
Bastian Harrach and Houcine Meftahi.
\newblock Global uniqueness and {L}ipschitz-stability for the inverse {R}obin
  transmission problem.
\newblock {\em SIAM Journal on Applied Mathematics}, 79(2):525--550, 2019.

\bibitem{YangSunChu20}
Liangliang Yang, Xiaogang Sun, and Yuanli Chu.
\newblock Boundary shape inversion of two-dimensional steady-state heat
  transfer system based on finite volume method and decentralized fuzzy
  adaptive {PID} control.
\newblock {\em Applied Sciences}, 10(1):153, 2020.

\bibitem{BucurGiacomini16}
Dorin Bucur and Alessandro Giacomini.
\newblock Shape optimization problems with robin conditions on the free
  boundary.
\newblock In {\em Annales de l'Institut Henri Poincar{\'e} C, Analyse non
  lin{\'e}aire}, volume~33, pages 1539--1568. Elsevier, 2016.

\bibitem{KarageorghisBin-MohsinLesnicEtAl15}
Andreas Karageorghis, B~Bin-Mohsin, D~Lesnic, and L~Marin.
\newblock Simultaneous numerical determination of a corroded boundary and its
  admittance.
\newblock {\em Inverse Problems in Science and Engineering}, 23(7):1120--1137,
  2015.

\bibitem{Rundell08}
William Rundell.
\newblock Recovering an obstacle and its impedance from cauchy data.
\newblock {\em Inverse problems}, 24(4):045003, 2008.

\bibitem{CakoniKressSchuft10a}
Fioralba Cakoni, Rainer Kress, and Christian Schuft.
\newblock Simultaneous reconstruction of shape and impedance in corrosion
  detection.
\newblock {\em Methods and Applications of Analysis}, 17(4):357--378, 2010.

\bibitem{CakoniKressSchuft10b}
Fioralba Cakoni, Rainer Kress, and Christian Schuft.
\newblock Integral equations for shape and impedance reconstruction in
  corrosion detection.
\newblock {\em Inverse Problems}, 26(9):095012, 2010.

\bibitem{CakoniKress12}
Fioralba Cakoni and Rainer Kress.
\newblock Integral equation methods for the inverse obstacle problem with
  generalized impedance boundary condition.
\newblock {\em Inverse Problems}, 29(1):015005, 2012.

\bibitem{CakoniHuKress14}
Fioralba Cakoni, Yuqing Hu, and Rainer Kress.
\newblock Simultaneous reconstruction of shape and generalized impedance
  functions in electrostatic imaging.
\newblock {\em Inverse Problems}, 30(10):105009, 2014.

\bibitem{KaipioSomersalo05}
Jari Kaipio and Erkki Somersalo.
\newblock {\em Statistical and Computational Inverse Problems}, volume 160 of
  {\em Applied Mathematical Sciences}.
\newblock Springer-Verlag New York, 2005.

\bibitem{CalvettiSomersalo07}
Daniela Calvetti and Erkki Somersalo.
\newblock {\em An introduction to {B}ayesian scientific computing: {T}en
  lectures on subjective computing}, volume~2.
\newblock Springer Science \& Business Media, 2007.

\bibitem{Stuart10}
Andrew~M. Stuart.
\newblock Inverse problems: {A B}ayesian perspective.
\newblock {\em Acta Numerica}, 19:451--559, 2010.

\bibitem{Wong01}
R.~Wong.
\newblock {\em Asymptotic Approximations of Integrals}.
\newblock Society for Industrial and Applied Mathematics, 2001.

\bibitem{NocedalWright06}
Jorge Nocedal and Stephen Wright.
\newblock {\em Numerical optimization}.
\newblock Springer Science \& Business Media, 2006.

\bibitem{DardeHyvonenSeppanenEtAl13a}
J{\'e}r{\'e}mi Dard{\'e}, Nuutti Hyv{\"o}nen, Aku Sepp{\"a}nen, and Stratos
  Staboulis.
\newblock Simultaneous reconstruction of outer boundary shape and admittivity
  distribution in electrical impedance tomography.
\newblock {\em SIAM Journal on Imaging Sciences}, 6(1):176--198, 2013.

\bibitem{DelfourZolesio11}
Michel~C Delfour and J-P Zol{\'e}sio.
\newblock {\em Shapes and geometries: metrics, analysis, differential calculus,
  and optimization}.
\newblock SIAM, 2011.

\bibitem{KolehmainenLassasOla05}
Ville Kolehmainen, Matti Lassas, and Petri Ola.
\newblock The inverse conductivity problem with an imperfectly known boundary.
\newblock {\em SIAM Journal on Applied Mathematics}, 66(2):365--383, 2005.

\bibitem{KolehmainenLassasOlaEtAl13}
Ville Kolehmainen, Matti Lassas, Petri Ola, and Samuli Siltanen.
\newblock Recovering boundary shape and conductivity in electrical impedance
  tomography.
\newblock {\em Inverse Problems \& Imaging}, 7(1):217, 2013.

\bibitem{KolehmainenLassasOla07}
Ville Kolehmainen, Matti Lassas, and Petri Ola.
\newblock Calder{\'o}n’s inverse problem with an imperfectly known boundary
  in two and three dimensions.
\newblock In {\em J. Phys.: Conf. Ser}, volume~73, page 012002, 2007.

\bibitem{BorziSchulz11}
Alfio Borz{\`\i} and Volker Schulz.
\newblock {\em Computational optimization of systems governed by partial
  differential equations}.
\newblock SIAM, 2011.

\bibitem{HinzePinnauUlbrichEtAl08}
Michael Hinze, Ren{\'e} Pinnau, Michael Ulbrich, and Stefan Ulbrich.
\newblock {\em Optimization with {PDE} constraints}, volume~23.
\newblock Springer Science \& Business Media, 2008.

\bibitem{Troltzsch10}
Fredi Tr{\"o}ltzsch.
\newblock {\em Optimal control of partial differential equations: {T}heory,
  methods, and applications}, volume 112.
\newblock American Mathematical Soc., 2010.

\bibitem{BoyleAdlerLionheart12}
Alistair Boyle, Andy Adler, and William~RB Lionheart.
\newblock Shape deformation in two-dimensional electrical impedance tomography.
\newblock {\em IEEE transactions on medical imaging}, 31(12):2185--2193, 2012.

\bibitem{DardeHyvonenSeppanenEtAl13b}
J{\'e}r{\'e}mi Dard{\'e}, Nuutti Hyv{\"o}nen, Aku Sepp{\"a}nen, and Stratos
  Staboulis.
\newblock Simultaneous recovery of admittivity and body shape in electrical
  impedance tomography: An experimental evaluation.
\newblock {\em Inverse Problems}, 29(8):085004, 2013.

\bibitem{HyvonenKarhunenSeppanen10}
Nuutti Hyv{\"o}nen, Kimmo Karhunen, and Aku Sepp{\"a}nen.
\newblock Fr{\'e}chet derivative with respect to the shape of an internal
  electrode in electrical impedance tomography.
\newblock {\em SIAM Journal on Applied Mathematics}, 70(6):1878--1898, 2010.

\bibitem{Sylvester90}
John Sylvester.
\newblock An anisotropic inverse boundary value problem.
\newblock {\em Communications on Pure and Applied Mathematics}, 43(2):201--232,
  1990.

\bibitem{Isakov06}
Victor Isakov.
\newblock {\em Inverse problems for partial differential equations}, volume
  127.
\newblock Springer, 2006.

\bibitem{LassasSaksmanSiltanen09}
Matti Lassas, Eero Saksman, and Samuli Siltanen.
\newblock Discretization-invariant bayesian inversion and besov space priors.
\newblock {\em Inverse Problems \& Imaging}, 3(1):87, 2009.

\bibitem{Bui-ThanhGhattasMartinEtAl13}
Tan Bui-Thanh, Omar Ghattas, James Martin, and Georg Stadler.
\newblock A computational framework for infinite-dimensional {B}ayesian inverse
  problems {P}art {I}: {T}he linearized case, with application to global
  seismic inversion.
\newblock {\em SIAM Journal on Scientific Computing}, 35(6):A2494--A2523, 2013.

\bibitem{DelfourZolesio01}
M.~C. Delfour and J.~P. Zol\'{e}sio.
\newblock {\em Shapes and geometries: analysis, differential calculus, and
  optimization}.
\newblock Society for Industrial and Applied Mathematics, 2001.

\bibitem{Kato13}
Tosio Kato.
\newblock {\em Perturbation theory for linear operators}, volume 132.
\newblock Springer Science \& Business Media, 2013.

\bibitem{LechleiterRieder08}
Armin Lechleiter and Andreas Rieder.
\newblock Newton regularizations for impedance tomography: {C}onvergence by
  local injectivity.
\newblock {\em Inverse problems}, 24(6):065009, 2008.

\bibitem{Daubechies92}
Ingrid Daubechies.
\newblock {\em Ten lectures on wavelets}.
\newblock SIAM, 1992.

\bibitem{FlathWilcoxAkcelikEtAl11}
Pearl~H. Flath, Lucas~C. Wilcox, Volkan Ak\c{c}elik, Judy Hill, Bart van
  Bloemen~Waanders, and Omar Ghattas.
\newblock Fast algorithms for {B}ayesian uncertainty quantification in
  large-scale linear inverse problems based on low-rank partial {H}essian
  approximations.
\newblock {\em SIAM Journal on Scientific Computing}, 33(1):407--432, 2011.

\bibitem{PetraMartinStadlerEtAl14}
Noemi Petra, James Martin, Georg Stadler, and Omar Ghattas.
\newblock A computational framework for infinite-dimensional {B}ayesian inverse
  problems: {P}art {II}. {S}tochastic {N}ewton {MCMC} with application to ice
  sheet inverse problems.
\newblock {\em SIAM Journal on Scientific Computing}, 36(4):A1525--A1555, 2014.

\bibitem{KhristenkoScarabosioSwierczynskiEtAl19}
U~Khristenko, L~Scarabosio, P~Swierczynski, E~Ullmann, and B~Wohlmuth.
\newblock Analysis of boundary effects on {PDE}-based sampling of
  {W}hittle--{M}at\'{e}rn random fields.
\newblock {\em SIAM/ASA Journal on Uncertainty Quantification}, 7(3):948--974,
  2019.

\bibitem{DaonStadler18}
Yair Daon and Georg Stadler.
\newblock Mitigating the influence of boundary conditions on covariance
  operators derived from elliptic {PDEs}.
\newblock {\em Inverse Problems and Imaging}, 12(5):1083--1102, 2018.

\bibitem{RoininenHuttunenLasanen14}
Lassi Roininen, Janne M.~J. Huttunen, and Sari Lasanen.
\newblock Whittle-mat\'{e}rn priors for {B}ayesian statistical inversion with
  applications in electrical impedance tomography.
\newblock {\em Inverse Problems \& Imaging}, 8(2):561, 2014.

\bibitem{VillaPetraGhattas20}
Umberto {Villa}, Noemi {Petra}, and Omar {Ghattas}.
\newblock {hIPPYlib: An Extensible Software Framework for Large-Scale Inverse
  Problems Governed by PDEs; Part I: Deterministic Inversion and Linearized
  Bayesian Inference}.
\newblock {\em arXiv e-prints}, 2020.
\newblock In review.

\bibitem{KovalAlexanderianStadler20}
Karina Koval, Alen Alexanderian, and Georg Stadler.
\newblock Optimal experimental design under irreducible uncertainty for linear
  inverse problems governed by pdes.
\newblock {\em Inverse Problems}, 2020.

\bibitem{GirolamiCalderhead11}
Mark Girolami and Ben Calderhead.
\newblock Riemann manifold {L}angevin and {H}amiltonian {M}onte {C}arlo
  methods.
\newblock {\em Journal of the Royal Statistical Society: Series B (Statistical
  Methodology)}, 73(2):123--214, 2011.

\bibitem{BuiGirolami14}
Tan Bui-Thanh and Mark Girolami.
\newblock Solving large-scale pde-constrained bayesian inverse problems with
  riemann manifold hamiltonian monte carlo.
\newblock {\em Inverse Problems}, 30(11):114014, 2014.

\bibitem{RobertsTweedie96}
Gareth~O Roberts and Richard~L Tweedie.
\newblock Exponential convergence of {L}angevin distributions and their
  discrete approximations.
\newblock {\em Bernoulli}, 2(4):341--363, 1996.

\bibitem{roberts1998optimal}
Gareth~O Roberts and Jeffrey~S Rosenthal.
\newblock Optimal scaling of discrete approximations to {L}angevin diffusions.
\newblock {\em Journal of the Royal Statistical Society: Series B (Statistical
  Methodology)}, 60(1):255--268, 1998.

\bibitem{ApteHairerStuartEtAl07}
Amit Apte, Martin Hairer, AM~Stuart, and Jochen Voss.
\newblock Sampling the posterior: {A}n approach to non-{G}aussian data
  assimilation.
\newblock {\em Physica D: Nonlinear Phenomena}, 230(1-2):50--64, 2007.

\bibitem{atchade2006adaptive}
Yves~F. Atchadé.
\newblock An adaptive version for the {M}etropolis adjusted {L}angevin
  algorithm with a truncated drift.
\newblock {\em Methodology and Computing in applied Probability}, 8:235--254,
  2006.

\bibitem{HaarioSaksmanTamminen01}
Heikki Haario, Eero Saksman, Johanna Tamminen, et~al.
\newblock An adaptive {M}etropolis algorithm.
\newblock {\em Bernoulli}, 7(2):223--242, 2001.

\bibitem{AndrieuThoms08}
Christophe Andrieu and Johannes Thoms.
\newblock A tutorial on adaptive {MCMC}.
\newblock {\em Statistics and computing}, 18(4):343--373, 2008.

\bibitem{gelman1992inference}
Andrew Gelman, Donald~B Rubin, et~al.
\newblock Inference from iterative simulation using multiple sequences.
\newblock {\em Statistical science}, 7(4):457--472, 1992.

\bibitem{brooks1998general}
Stephen~P Brooks and Andrew Gelman.
\newblock General methods for monitoring convergence of iterative simulations.
\newblock {\em Journal of computational and graphical statistics},
  7(4):434--455, 1998.

\bibitem{flegal2008markov}
James~M Flegal, Murali Haran, and Galin~L Jones.
\newblock Markov chain {M}onte {C}arlo: Can we trust the third significant
  figure?
\newblock {\em Stat. Sci.}, pages 250--260, 2008.

\bibitem{vats2019multivariate}
Dootika Vats, James~M Flegal, and Galin~L Jones.
\newblock Multivariate output analysis for {M}arkov chain {M}onte {C}arlo.
\newblock {\em Biometrika}, 106(2):321--337, apr 2019.

\bibitem{NiskanenDazelGrobyEtAl19}
Matti Niskanen, Olivier Dazel, Jean-Philippe Groby, Aroune Duclos, and Timo
  L{\"a}hivaara.
\newblock Characterising poroelastic materials in the ultrasonic range - {A}
  {B}ayesian approach.
\newblock {\em Journal of Sound and Vibration}, 456:30--48, 2019.

\bibitem{KaipioSomersalo07}
Jari Kaipio and Erkki Somersalo.
\newblock Statistical inverse problems: {D}iscretization, model reduction and
  inverse crimes.
\newblock {\em Journal of Computational and Applied Mathematics},
  198(2):493--504, 2007.

\bibitem{Grafakos08}
Loukas Grafakos.
\newblock {\em Classical {F}ourier analysis}, volume~2.
\newblock Springer, 2008.

\bibitem{EisenstatWalker96}
Stanley~C Eisenstat and Homer~F Walker.
\newblock Choosing the forcing terms in an inexact {N}ewton method.
\newblock {\em SIAM Journal on Scientific Computing}, 17(1):16--32, 1996.

\end{thebibliography}

\end{document}